\documentclass[11pt]{amsart}
\usepackage{graphicx} 

\usepackage[english]{babel}
\usepackage{import,xspace}
\usepackage{amsmath,amsfonts,amsthm,amssymb,stmaryrd,mathtools}

\usepackage[alphabetic,abbrev]{amsrefs} 
\renewcommand{\MR}[1]{}

\usepackage{url}
\usepackage{fancyhdr}
\usepackage{graphicx}
\usepackage{microtype}
\usepackage{enumitem}
\usepackage{wrapfig, caption}
\usepackage[hidelinks]{hyperref}
\usepackage{xcolor}
\usepackage{geometry}
\usepackage{dutchcal}

\hypersetup{
    colorlinks,
    citecolor=black!40!blue,
    urlcolor=black!40!blue,
    linkcolor=black!40!blue,
}

\usepackage{tikz}
\usetikzlibrary{cd,calc}
\usetikzlibrary{arrows}
\usetikzlibrary{decorations.markings}
\newcommand{\pd}[2][black]{\filldraw[#1] (#2) circle (1.5pt)} 
\newcommand{\pdg}[3]{
    \filldraw[black] (#1) circle (1.5pt) node[anchor=#2,inner sep=0] {$\scriptstyle{\strut #3}$}
}

\newcommand\crossing{\tikz[baseline=(origin.base),scale=0.5]{
    \coordinate (origin) at (0,0); 
    \pd{0,0}; \pd{.5,0}; \pd{0,.5}; \pd{.5,.5};
    \draw (0,0) -- (.5,.5);
    \draw (.5,0) -- (0,.5);
}}
\newcommand\mult{\tikz[baseline=(origin.base),scale=0.5]{
    \coordinate (origin) at (0,0); 
    \pd{.25,0}; \pd{0,.5}; \pd{.5,.5};
    \draw (.25,0) -- (.25,.25) -- (0,.5);
    \draw (.25,.25) -- (.5,.5);
}}
\newcommand\comult{\tikz[baseline=(origin.base),scale=0.5]{
    \coordinate (origin) at (0,0); 
    \pd{.25,.5}; \pd{0,0}; \pd{.5,0};
    \draw (.25,.5) -- (.25,.25) -- (0,0);
    \draw (.25,.25) -- (.5,0);
}}
\newcommand\coarsener{\tikz[baseline=(origin.base),scale=0.5]{
    \coordinate (origin) at (0,0);
    \pd{0,0}; \pd{.5,0}; \pd{1,0}; \pd{0,.5}; \pd{.5,.5}; \pd{1,.5};
    \draw (0,0) -- (0,.5); \draw (.5,0) -- (.5,.5); \draw (1,0) -- (1,.5);
    \draw (0,.25) -- (.3,.25) to[out=up,in=up] (.7,.25) -- (1,.25);
}}
\newcommand\upperdots[1]{\tikz[baseline=(origin.base),scale=0.5]{
    \coordinate (origin) at (0,0);
    \foreach \i in {0,...,\numexpr#1-1\relax} {
        \pd{0.5*\i,.5};
        \draw (0.5*\i,.3) -- (0.5*\i,.5);
    }
}}
\newcommand\cuppartition{\tikz[baseline=(origin.base),scale=0.5]{
    \coordinate (origin) at (0,0); 
    \pd{0,.5}; \pd{.5,.5};
    \draw (0,.5) to[in=down,out=down] (.5,.5);
}}
\newcommand\verticaldots{\tikz[baseline=(origin.base),scale=0.5]{
    \coordinate (origin) at (0,0);
    \pd{0,0};\pd{0,.5};
    \draw (0,0) -- (0,.18); \draw (0,.32) -- (0,.5);
}}

\newcommand\psmat[1]{\left(\begin{smallmatrix}#1\end{smallmatrix}\right)}

\usepackage{zref-clever} 
\newcommand\Cref[1]{\zcref[S]{#1}} 
\newcommand\newtheoremx[2]{%
\AddToHook{env/#1/begin}{
\zcsetup{countertype={theorem=#1}}}
\newtheorem{#1}[theorem]{#2}
}

\newtheorem{theorem}{Theorem}[section]
\newtheoremx{claim}{Claim}
\newtheoremx{proposition}{Proposition}
\newtheoremx{corollary}{Corollary}
\newtheoremx{lemma}{Lemma}

\newtheorem{introtheorem}{Theorem}

\zcsetup{countertype={introtheorem=theorem}}

\theoremstyle{definition}
\newtheoremx{definition}{Definition}
\newtheoremx{notation}{Notation}

\zcRefTypeSetup{construction}{
  Name-sg = Construction,
  Name-pl = Constructions,
}
\newtheoremx{construction}{Construction}

\zcRefTypeSetup{assumption}{
  Name-sg = Assumption,
  Name-pl = Assumptions,
}
\newtheoremx{assumption}{Assumption}

\theoremstyle{remark}
\newtheoremx{example}{Example}
\newtheoremx{remark}{Remark}
\newtheoremx{question}{Question}

\title[Abelian envelopes and quantum groups in Knop's interpolation categories]{Abelian envelopes and quantum groups\\in Knop's interpolation categories}

\author{Johannes Flake}
\address{Mathematical Institute, University of Bonn, Endenicher Allee 60, 53115 Bonn, Germany}
\email{flake@math.uni-bonn.de}

\author{Peter Mader}
\address{Mathematical Institute, University of Bonn, Endenicher Allee 60, 53115 Bonn, Germany}
\email{mader@math.uni-bonn.de}

\date{}

\newcommand\cA{\mathcal A}
\newcommand\cB{\mathcal B}
\newcommand\cC{\mathcal C}
\newcommand\cD{\mathcal D}

\newcommand\cN{\mathcal N}
\newcommand\cT{\mathcal T}
\newcommand\cS{\mathcal S}
\newcommand\cV{\mathcal V}
\newcommand\fA{\mathfrak A}
\newcommand\fB{\mathfrak B}
\newcommand\kk{\Bbbk}
\newcommand\Rep{\operatorname{Rep}}

\newcommand\uRep{\underline{\smash{\operatorname{Rep}}}}
\newcommand\uPerm{\underline{\smash{\operatorname{Perm}}}}
\newcommand\Id{\operatorname{Id}}
\newcommand\id{\operatorname{id}}
\newcommand\im{\operatorname{im}}
\newcommand\Ind{\operatorname{Ind}}
\newcommand\Pro{\operatorname{Pro}}
\newcommand\Surj{\operatorname{Surj}}
\newcommand\Aut{\operatorname{Aut}}
\newcommand\Comod{\operatorname{Comod}}
\renewcommand\:\colon
\renewcommand\o\otimes
\newcommand\one{\mathbf{1}}

\newcommand\ol\overline
\newcommand\GL{\mathrm{GL}}
\newcommand\Aff{\mathrm{Aff}}
\newcommand\Fun{\operatorname{Fun}}

\newcommand\Hom{\operatorname{Hom}}

\renewcommand\Vec{\mathrm{Vec}}
\newcommand\cG{\mathcal G}
\newcommand\ZZ{\mathbb Z}
\newcommand\NN{\mathbb N}
\newcommand\CC{\mathbb C}
\newcommand\FF{\mathbb F}
\renewcommand\o\otimes
\newcommand\cO{\mathcal O}
\newcommand\ev{\operatorname{ev}}
\newcommand\coev{\operatorname{coev}}

\renewcommand\:\colon

\newcommand\Rel{\operatorname{Rel}}

\newcommand\FinSet{\mathrm{FinSet}}

\newcommand\op{^{\mathrm{op}}}
\newcommand\ab{^{\mathrm{ab}}}
\newcommand\kar{^{\mathrm{kar}}}

\newcommand\topp{\mathrm{top}}
\newcommand\bott{\mathrm{bot}}
\newcommand\eps{\varepsilon}
\newcommand\vfi{\varphi}
\newcommand\subobj{\mathrm{sub}}
\newcommand\End{\operatorname{End}}
\newcommand\abs[1]{\lvert{}#1\rvert{}}
\newcommand\sspan{\operatorname{span}}
\newcommand\angl[1]{\langle{}#1\rangle{}} 
\newcommand\gtcategory[1]{\angl{#1}_{\mathrm{gt}}}
\newcommand\tto{\twoheadrightarrow}

\begin{document}

\begin{abstract} We show the existence of abelian envelopes for large classes of subcategories of Knop's tensor envelopes. It follows that these abelian envelopes correspond to groups or more generally Hopf algebras (``quantum groups'') in their parent interpolation categories. As a special case, we obtain that the interpolation categories in the sense of Flake--Maassen (2022) of all so-called group-theoretical easy quantum groups have abelian envelopes given by quantum groups in Deligne's interpolation category of the symmetric groups. Similarly, we apply our theory to interpolation categories of wreath products, finite general linear groups, and finite affine groups. Our main technique is an explicit construction of objects corresponding to coordinate algebras using the calculus of relations that Knop's categories are based on.
\end{abstract}

\maketitle

\section{Introduction}

A \emph{tensor category} is an abelian category with a tensor product satisfying some natural axioms that apply, for instance, to categories of representations of finite or algebraic groups. In fact, these examples can be re-interpreted as an organizational principle: for any symmetric tensor category $\cC$, symmetric tensor functors to $\cC$ are the same thing as groups, in a suitable general sense, internal to $\cC$ by the theory of Tannakian reconstruction. For any group, the corresponding tensor category is given by the category of representations. More generally, a certain class of functors that are not necessarily symmetric corresponds to Hopf algebras internal to $\cC$ with their categories of comodules. In the spirit of \cite{Drinfeld}, we speak of these as quantum groups internal to $\cC$ and their representations.\footnote{Note that the term quantum group can mean various other things, too, such as quantized universal enveloping algebra.}

A deep theorem by Deligne \cite{Del-tensorielles} establishes that, over an algebraically closed field in characteristic zero, a large class of tensor categories is given by a group internal to the category of vector spaces or super vector spaces. This is the class of those symmetric tensor categories in which the (composition) lengths of tensor powers grow at most exponentially. Hence, to study tensor categories beyond the scope of Deligne's theorem, one has to turn to positive characteristics or to superexponential growth, and both directions have led to active areas of research.

Among the tensor categories of superexponential growth in characteristic zero, an important class is given by \emph{Knop's interpolation categories}, also known as \emph{Knop's tensor envelopes} \cite{knop}. These are tensor categories defined by linearizing the calculus of relations of a regular category $\cA$. For each $\cA$, a family of categories $\cT(\cA,\delta)$ is constructed, parametrized by so-called degree functions $\delta$ that control the linearizations. The categories $\cT(\cA,\delta)$ are generally non-abelian versions of tensor categories called pseudo-tensor categories. For suitable $\cA$, they are semisimple tensor categories for almost all (generic) choices of $\delta$, and admit semisimple quotient tensor categories for the remaining (special) choices of $\delta$. As the semisimple quotients can be identified with categories of representations of finite groups, one can view the $\cT(\cA,\delta)$ as interpolation categories for the respective family of finite groups. For instance, for $\cA$ the opposite category of finite sets, Deligne's interpolation categories $\uRep(S_t)$ for the symmetric groups \cite{deligne2007repst} are recovered by Knop's construction.

Note that besides recovering $\uRep(S_t)$, Knop's interpolation categories are also closely related to Harman--Snowden's tensor categories constructed from pro-oligomorphic groups \cite{oligotensor} by work of Snowden \cite{knoptooligo}, as we discuss in \Cref{ss:oligo}.

Instead of passing to the semisimple quotient at the special parameters $\delta$, one can also ask if the categories $\cT(\cA,\delta)$ define a universal tensor category in which they embed. This leads to the notion of \emph{abelian envelopes} \cites{Cou-abelian,BEO,CEOP,flake2026monoidaladjunctionsabelianenvelopes}, which has become a central concept in the theory of tensor categories. Deciding if such an abelian envelope exists for a given pseudo-tensor category, and then constructing it, is however often a difficult problem.

In this paper, we solve this problem for a large class of suitable subcategories of Knop's interpolation categories. Using the calculus of relations, we show that these categories have abelian envelopes which are given by (quantum) groups in the parent interpolation category by a suitable version of Tannakian reconstruction.

\medskip

More concretely, fix an algebraically closed base field of characteristic $0$ and let $\Rel(\cA)$ be the monoidal category of relations of a fixed regular category $\cA$ subject to some assumptions, with tensor products of morphisms given by products and composition given by pullbacks followed by images. Then Knop's tensor envelope $\cT(\cA,\delta)$ is a linearization of $\Rel(\cA)$ specified by the parameter $\delta$.
We define in \Cref{s:group-theoretical-subcats} a class of \emph{group-theoretical} subcategories of $\Rel(\cA)$. Each such subcategory $\cC$ generates a non-full pseudo-tensor subcategory $\cC_\delta$ of $\cT(\cA,\delta)$ with an inclusion functor $\iota$.

We first show that $\cC_\delta$ and $\cT_\delta:=\cT(\cA,\delta)$ have comparable sets of special parameters $\delta$.

\begin{introtheorem}[\Cref{t:c-delta-semisimple}] \label{t:intro-semisimple} 
    $\cC_\delta$ is semisimple whenever $\cT_\delta$ is, and the opposite implication holds if $\cA$ is Barr-exact and Malcev.
\end{introtheorem}

Then we examine conditions under which an abelian envelope of the Knop category $\cT_\delta$ gives rise to an abelian envelope of the subcategory $\cC_\delta$. For an injection $f: X \hookrightarrow Y$ in $\cA$, a \emph{$\cC$-rank object} is an object $R \in \cA$ such that there is a relation of the form $Y \xhookleftarrow{f} X \to R$ in $\cC$ and for all relations $Y \xhookleftarrow{f} X \stackrel{\eps}{\twoheadleftarrow} W \to R$ in $\cC$, the morphism $\eps$ is an isomorphism.

\begin{introtheorem}[\Cref{prop:adjoint}, \Cref{t:group-theoretical-ab-env}, \Cref{t:group-theoretical-comod}] \label{t:intro-adjoint}
    Assume that the base field has characteristic $0$ and every injection in $\cA$ has a $\cC$-rank object. Then the functor
    \[
        \Ind(\iota): \Ind(\cC_\delta) \longrightarrow \Ind(\cT_\delta)
    \]
    has a right adjoint. Moreover, if $\cT_\delta$ has a global splitting object of the form $\iota(S)$ with $S \in \cC_\delta$, then $\cC_\delta$ and $\cT_\delta$ admit abelian envelopes $\cC_\delta\ab$, $\cT_\delta\ab$ with enough projectives, and $\cC_\delta\ab$ is equivalent to a category of representations of a (quantum) group internal to $\cT_\delta\ab$.
\end{introtheorem}

Here, a \emph{global splitting object} is a non-zero object tensoring with which turns any morphism into a split morphism. Splitting objects play a central role for constructing abelian envelopes. By an observation from \cite{flake2026monoidaladjunctionsabelianenvelopes}, they can be transferred along certain adjunctions (\Cref{l:adjoint-preserves-splitting-objects}). This allows us to deduce the existence of the abelian envelope of $\cC_\delta$ from the existence of the asserted adjoint functor. Its existence, in turn, is proved by an explicit construction of an object that is the image of the tensor unit under the adjoint (\Cref{p:representing-object}). This construction boils down to a careful analysis of the calculus of relations.

It follows that all asserted abelian envelopes exist and by their universal properties, we obtain an (exact) tensor functor $\cC_\delta\ab \to \cT_\delta\ab$ extending $\iota$. By a suitable version of Tannakian reconstruction (\Cref{thm:rec-1}), this tensor functor corresponds to a quantum group internal to $\cT_\delta\ab$, or even a group, if $\cC$ is symmetric as a monoidal subcategory of $\Rel(\cA)$. We discuss an explicit description of these quantum groups in \Cref{ss:quantum-groups}.

\medskip

We apply \Cref{t:intro-adjoint} to several Knop categories, with Deligne's category $\uRep(S_t)$ for $t \in \CC$ as the main example. In this case, the category $\Rel(\cA)$ is based on the calculus of partition diagrams. Some of its subcategories have been studied and classified in the context of so-called ``easy quantum groups'' (see e.g.~\cites{banica2009liberation,raum2016full}). The terminology \emph{group-theoretical} comes from this classification. The categories $\cC_\delta$ correspond to the categories $\uRep(\cC, t)$ from \cite{flake2022semisimplicity}, where the authors parametrize their indecomposable objects and give a semisimplicity criterion which for group-theoretical $\cC$ agrees with \Cref{t:intro-semisimple}.

The category $\uRep(S_t)$ has a global splitting object $\Delta$ which is contained in $\uRep(\cC, t)$ for each of the uncountably many group-theoretical categories of partitions $\cC$. We can apply \Cref{t:intro-adjoint} and \Cref{thm:one-atom} to obtain:

\begin{introtheorem}[\Cref{c:partitions-ab-env}]
    For every group-theoretical category $\cC$ of partitions, the pseudo-tensor category $\uRep(\cC, t)$ has an abelian envelope with enough projectives, given by the category of representations of a (quantum) group internal to $\uRep(S_t)\ab$.
\end{introtheorem}

As a specific example, we consider the subcategory $\cC$ given by partitions with even block sizes. It is symmetric, and we show that it corresponds to the group scheme $(S_2)^t$ in $\uRep(S_t)\ab$, see \Cref{p:hopf-alg-iso}. Note that $t$ does not need to be an integer here, and the non-integral power of $S_2$ is interpreted using Etingof's non-integral tensor powers \cite{reptheocomplexrank}*{\S\S~3.2 and~4.1}.

While the calculus of relations $\Rel(\cA)$ only depends on $\cA$ up to equivalence of categories, the class of its group-theoretical subcategories can change even when $\cA$ is replaced by an equivalent category. We briefly explore in \Cref{ss:other-finset-op} the consequences of replacing the calculus of partitions in $\uRep(S_t)$ by a different but equivalent calculus; in particular, we observe that one obtains more examples of groups in $\uRep(S_t)\ab$ in this way (\Cref{ex:klein-four-group}).

Knop's construction also produces interpolation categories for wreath products $(G \wr S_n)_{n\ge0}$, where $G$ is a fixed finite group, for general linear groups $(\GL_n(\FF_q))_{n\ge0}$, or for finite affine groups $(\Aff_n(\FF_q))_{n\ge0}$, where $\FF_q$ is a finite field. They are obtained by choosing as the regular category $\cA$ the opposite category of $G\text{-}\mathrm{FreeFinSet}$ of finite free $G$-sets, the category of finite-dimensional vector spaces $\mathrm{Vect}_{\FF_q}^{\mathrm{fd}}$, or the category of non-empty finite-dimensional affine spaces $\mathrm{Aff}_{\FF_q}^{\mathrm{fd}}$; to be more precise, we choose natural skeletal versions of these categories as $\cA$. We show:

\begin{introtheorem}[\Cref{c:wreath-products-ab-env}, \Cref{p:linear-affine-ab-env}]
    For $\cA$ any of $G\text{-}\mathrm{FreeFinSet}\op$, $\mathrm{Vect}_{\FF_q}^{\mathrm{fd}}$, $\mathrm{Aff}_{\FF_q}^{\mathrm{fd}}$, with $q>2$ in the latter two cases, and for all group-theoretical subcategories $\cC$ of $\Rel(\cA)$, the corresponding pseudo-tensor subcategories $\cC_\delta \hookrightarrow \cT_\delta=\cT(\cA,\delta)$ have abelian envelopes with enough projectives given by categories of representations of (quantum) groups internal to $\cT_\delta\ab$.
\end{introtheorem}

Besides establishing the existence of a host of new and interesting abelian envelopes using our novel uniform criteria, including for the uncountably many group-theoretical subcategories of $\uRep(S_t)$, we think that our work can lay the foundation for a theory of (quantum) groups internal to interpolation categories.
Note that Harman--Kalinov's classification of simple associative algebras in $\uRep(S_t)$ for $t\notin\NN_0$ in \cite{HK}, which we regard as another contribution towards such a theory, does, for instance, not cover the coordinate algebra of our example $(S_2)^t$, as it has infinite length.

\subsection*{Use of AI} We started working on this project in December 2025. Almost all results of this paper were obtained without using AI, and the paper was written by hand. We did use {GPT-6} Astra for the exact statements and proofs of the probably known statements in \Cref{ss:tannaka-duality}, which we could not find in the literature in this exact form; to find the proofs of \Cref{p:subgroup-of-G} and \Cref{p:hopf-alg-iso}; and for proof-reading the final document. In particular, it helped us identify and fix a subtle issue with \Cref{d:group-theoretical}.

\subsection*{Acknowledgments} J.F.~thanks David Green, Nate Harman, Thorsten Heidersdorf, Alexander Mang, Andrew Snowden, and Moritz Weber for useful discussion, some of them from a couple of years ago. In particular, J.F. thanks Thorsten Heidersdorf and Andrew Snowden for pointing out the content of \Cref{p:knoptooligo}(a). J.F.~was supported by the Hausdorff Center for Mathematics Bonn. P.M.~was supported by the Bonn International Graduate Schools (BIGS) and the German Academic Scholarship Foundation (Studienstiftung des deutschen Volkes).

\tableofcontents

\section{Background}

We fix a field $\kk$.

\subsection{Partial meet-semilattices and the Möbius basis}

In this subsection, we repeat basic definitions and facts about posets with meets given in \cite{knop}*{\S~7}. In particular, we will need the following notions.

\begin{definition}
    Let $L$ be a poset.

    \begin{enumerate}
        \item[(a)] We call $L$ a \emph{meet-semilattice} if for all $u, v \in L$, the meet $u \wedge v$ exists.
        \item[(b)] Consider the extension $L^\varnothing \coloneqq L \sqcup \{\varnothing\}$ by a new minimum $\varnothing$, i.e.~we put $\varnothing < u$ for all $u \in L$. We call $L$ a \emph{partial meet-semilattice} if $L^\varnothing$ is a meet-semilattice.
        \item[(c)] If $L$ is a partial meet-semilattice, its \emph{Möbius algebra} $A(L)$ is the free $\kk$-vector space with basis $L$. Multiplication of basis elements $u, v \in L$ is given by $u \wedge v$ if the meet exists, and 0 otherwise.
        \item[(d)] A subset $M\subseteq L$ is called \emph{order ideal} in $L$ if for all $u\in L$, $v \in M$, $u \leq v$ implies $u \in M$.
    \end{enumerate}
\end{definition}

Consider a finite partial meet-semilattice $L$ and its Möbius function $\mu_L: L \times L \to \kk$ uniquely determined by $\mu_L(u,v)=0$ for $u\not\leq v$ and $\sum_{u \leq w \leq v} \mu_L(u, w) = \delta_{u,v}$ for $u\leq v$ in $L$. It gives rise to elements $x(v) \coloneqq \sum_{u \leq v} \mu_L(u, v) u \in A(L)$ for $v \in L$. They form a basis of $A(L)$ called the \emph{Möbius basis}, which satisfies the absorption property
\begin{equation} \label{e:moebius-basis-properties}
    x(u) v = \begin{cases}
        x(u) & \text{if } u \leq v, \\
        0 & \text{otherwise},
    \end{cases} \qquad\text{and}\qquad
    x(u) x(v) = \begin{cases}
        x(u) & \text{if } u = v, \\
        0 & \text{otherwise}.
    \end{cases}
\end{equation}
The following statement is a direct consequence of the definitions:

\begin{lemma} \label{l:moebius-function-restricts}
    Let $L$ be a finite partial meet-semilattice and $M \subseteq L$ be an order ideal in $L$. Then $M$ is a partial meet-semilattice and the Möbius functions $\mu_L, \mu_M$ coincide on $M$.
\end{lemma}

\begin{lemma}[see \cite{knop}*{Lemma~7.1}] \label{l:lindstroems-lemma}
    Consider a finite partial meet-semilattice $L$ and a function $\vfi: L \to \kk$, which we extend linearly to a map $\vfi: A(L) \to \kk$. Then
    \[
        \det((\vfi(r\wedge s))_{r, s \in L}) = \prod_{r \in L} \vfi(x(r)),
    \] 
    with the convention that $\vfi(r\wedge s):=0$ whenever the meet does not exist in $L$. 
\end{lemma}

\subsection{Tensor categories} \label{bg:tensorcat} For us, a \emph{pseudo-tensor category} is a hom-finite $\kk$-linear additive and idempotent-complete rigid monoidal category with $\End(\one)\cong\kk$. A \emph{tensor category} is a finite-length abelian pseudo-tensor category, which recovers the notion from \cite{EGNO}. For any pseudo-tensor category $\cC$, an \emph{abelian envelope} is a fully faithful embedding $\cC\hookrightarrow\cT$ into a tensor category $\cT$ such that, for any tensor category $\cT'$, restriction along this embedding induces an equivalence between the categories
$$
\Fun^{\text{ex}}(\cT,\cT') \to \Fun^{\text{faith}}(\cC,\cT')
$$
of exact / faithful linear monoidal functors, respectively.

If every non-zero morphism between indecomposable objects of a pseudo-tensor category $\cC$ is an isomorphism, then $\cC$ is a semisimple tensor category and its own abelian envelope.

If $\cC$ is a spherical pseudo-tensor category, we have a notion of a \emph{trace} of an endomorphism. A morphism $f: X \to Y$ in $\cC$ is called \emph{negligible} if $\mathrm{tr}(fg) = 0$ for all morphisms $g: Y \to X$ in $\cC$. The negligible morphisms form a tensor ideal.

\begin{proposition} \label{p:semisimplification}
    Assume that $\kk$ is algebraically closed and $\cC$ is a spherical pseudo-tensor category over $\kk$. Then $\cC$ is semisimple if and only if all nilpotent endomorphisms have trace 0 and the ideal of negligible morphisms vanishes.
\end{proposition}
\begin{proof}
    If $\cC$ is semisimple, it is in particular abelian, so nilpotent endomorphisms have trace~0. The negligible morphisms vanish by \cite{ssnegligibles}*{Corollaire~7.1.7}. The converse statement is proved in \cite{onsemisimplification}*{Theorem~2.6}.
\end{proof}

\subsection{Knop's tensor envelopes} \label{bg:knop}

In the following, we assume that $\cA$ is a subobject-finite regular category and that the terminal object $\one$ has no proper subobjects.\footnote{Here, we use Knop's convention for regularity \cite{knop}*{\S~2}: adjoining a new absolutely initial object $\varnothing$ yields a finitely complete regular category in the usual sense. Non-existing fiber products in $\cA$ are given by $\varnothing$. Empty relations correspond to zero morphisms in the linearization.} Monomorphisms are called \emph{injective}, and extremal epimorphisms are called \emph{surjective}.

Let $\delta$ be a degree function on $\cA$ as in \cite{knop}*{3.1.~Definition}, that is, a rule assigning a value in a fixed field $\kk$ to any surjective morphism in $\cA$ such that isomorphisms are assigned the value $1$, $\delta$ is multiplicative with respect to compositions, and invariant under pull-backs.

In \cite{knop}*{\S 3} Knop constructs a pseudo-tensor category $\cT(\cA,\delta)$ in the sense of \Cref{bg:tensorcat} for any $(\cA,\delta)$ as above. Each object $X \in \cA$ gives rise to an object in $\cT(\cA, \delta)$ denoted by $[X]$. The terminal object of $\cA$ and the unit object in $\cT(\cA, \delta)$ are both denoted by $\one$. It is shown:

\begin{theorem}[\cite{knop}*{\S~8}] \label{thm:knop} Assume $\kk$ is of characteristic $0$. If $\cT(\cA,\delta)$ is semisimple, then
    \[
        \sum_{X \in \subobj(Y)} \mu(X, Y) \delta(X \tto \one) \ne 0
    \]
    for all $Y \in \cA$, where $\mu$ is the Möbius function of the poset $\subobj(Y)$. The opposite implication holds if $\cA$ is Barr-exact and Malcev.
\end{theorem}

\subsection{Pseudo-diagrammatic categories} \label{ss:pseudo-diagrammatic-categories}

In the following, let $F: \cC \to \cD$ be a linear monoidal functor between linear monoidal categories. Write $\tilde{F}: \tilde{\cC} \to \tilde{\cD}$ for the $\Ind$-completion.

\begin{definition}
    The functor $F$ is called \emph{dominant} if every object $Y \in \cD$ is a direct summand of $F(X)$ for some $X \in \cC$.
\end{definition}

\begin{lemma}[\cite{flake2026monoidaladjunctionsabelianenvelopes}*{Lemma~5.7}] \label{l:ind-adjoint}
    Assume that $\cC, \cD$ are rigid and $F$ is dominant. Assume that there is an object $H \in \tilde{\cC}$ such that
    \[
        \tilde{\cD}(\tilde{F}(-), \one)|_{\cC\op} \cong \tilde{\cC}(-, H)|_{\cC\op}
    \]
    as functors $\cC\op \to \Vec$. Then $\tilde{F}$ has a right adjoint sending $\one \mapsto H$.
\end{lemma}

Recall that a morphism $f$ is called \emph{split} if there is an (anti-parallel) morphism $g$ such that $fgf=f$. A non-zero object $X$ is said to \emph{split} a morphism $f$ if $X\o f$ and $f\o X $ are split. A \emph{global splitting object} is an object that splits all morphisms in a given category.

\begin{lemma} \label{l:adjoint-preserves-splitting-objects}
    Assume that $\cC, \cD$ are rigid and that $\tilde{F}$ has a right adjoint $G$ such that $G(\one_\cD)$ has a non-zero direct summand $U \in \cC$. Let $f$ be a morphism in $\cC$ and $S \in \cC$ be an object with $S\o U\ne0$. If $F(S)$ splits $F(f)$, then $S \o U$ splits $f$.
\end{lemma}
\begin{proof}
    We view $\cC \subseteq \tilde{\cC}, \cD \subseteq \tilde{\cD}$ as full subcategories and abbreviate $H \coloneqq G(\one_\cD)$. We have isomorphisms
    \[
        \tilde{\cC}(X, G(F(Y))) \cong \tilde{\cD}(\tilde{F}(X), F(Y)) \cong \tilde{\cD}(\tilde{F}(Y^* \o X), \one_\cD) \cong \tilde{\cC}(Y^* \o X, H) \cong \tilde{\cC}(X, Y \o H)
    \]
    natural in $X \in \tilde{\cC}$ and $Y \in \cC$ and thus obtain an isomorphism of functors $G(F(-)) = - \o H$.
    
    Now consider a morphism $f: A \to B$ in $\cC$, an object $S \in \cC$, and a splitting $g: F(B) \otimes F(S) \to F(A) \o F(S)$ of $F(f)\o F(S)$ in $\cD$ and thus in $\tilde{\cD}$. Then
    \[
        G(g): B \o S \o H = G(F(B \o S)) \longrightarrow G(F(A \o S)) = A \o S \o H
    \]
    is a splitting of $G(F(f\o S)) = f \o S \o H$. Then the co-/restriction of $G(g)$ to $B \o S \o U \to A \o S \o U$ is a splitting of $f\o S\o U$ in $\cC$. Similarly, using both left and right duals we get natural isomorphisms $G(F(S)\o F(A))\cong S\o H\o A$. Applying $G$ to a splitting of $F(S)\o F(f)$ and restricting to the direct summand $S\o U$ gives a splitting of $S\o U\o f$.
\end{proof}

In \cite{flake2026monoidaladjunctionsabelianenvelopes}*{\S 3}, a class of $\oplus$-pseudo-diagrammatic pseudo-tensor categories is defined. It is shown in \cite{flake2026monoidaladjunctionsabelianenvelopes}*{\S 4.2} that, in particular, Knop's tensor envelopes as recalled in \Cref{bg:knop} are $\oplus$-pseudo-diagrammatic. This has the following consequences for us:


\begin{theorem} \label{thm:adjoint-ab-env}
    Assume that $F\:\cC\to\cD$ is a linear monoidal functor between pseudo-tensor categories such that

    \begin{enumerate}
        \item[(a)] $\cD$ is an $\oplus$-pseudo-diagrammatic pseudo-tensor category,
        \item[(b)] $\cD$ has a global splitting object of the form $F(S)$ for some $S \in \cC$,
        \item[(c)] $\tilde{F}$ has a right adjoint $G$, and
        \item[(d)] $\one$ is a direct summand in $G(\one_\cD)$.
    \end{enumerate}

    Then $S$ is a global splitting object in $\cC$, $\cC$ and $\cD$ admit abelian envelopes with enough projectives, and $F$ induces a tensor functor $F\ab\:\cC\ab\to\cD\ab$.
\end{theorem}

\begin{proof}
The existence of the splitting object for $\cC$ follows from \Cref{l:adjoint-preserves-splitting-objects} using (b), (c), and (d).

Write $H=G(\one_\cD)$ and let $\eta$ be the adjunction unit. By adjunction, $\tilde{\cC}(\one,H)\cong\cD(\one_\cD,\one_\cD)=\kk$, with $\eta_\one$ corresponding to the identity. Thus (d) implies that $\eta_\one$ is a split monomorphism. The projection formula $GF(X)\cong X\o H$ identifies $\eta_X$ with $X\o\eta_\one$ for all $X\in\cC$, so the unit splits naturally on $\cC$. This implies that $F$ is faithful.

Consider the following kind of exactness condition for a pseudo-tensor category: any non-zero morphism $u\:U\to\one$ is a cokernel of the morphism $(u\o U-U\o u)$. This condition holds for $\cD$ by \cite{flake2026monoidaladjunctionsabelianenvelopes}*{Proposition~3.10} using (a). We claim it also holds in $\cC$.

Indeed, the above condition is equivalent to exactness of the image under $\Hom(-,X)$, for every $X$, of the sequence
$$
\Sigma_u \: U\o U\xrightarrow{u\o U-U\o u} U\xrightarrow{u} \one\to 0 .
$$
The splitting of the unit makes $\Hom(\Sigma_u,X)$ a retract of $\Hom(\Sigma_u,GF(X))$. By adjunction, the latter is isomorphic to $\Hom(F(\Sigma_u),F(X))$, which is exact by the corresponding condition for $\cD$, since $F(u)\ne0$. Hence $\Hom(\Sigma_u,X)$ is exact as well.

Now \cite{flake2026monoidaladjunctionsabelianenvelopes}*{Theorem~6.1} gives the claimed abelian envelopes and the universal property of abelian envelopes gives the asserted functor between them.
\end{proof}


\subsection{Tannakian reconstruction} \label{ss:tannaka-duality}

We start by recalling some very general, and probably known, ideas related to Tannakian reconstruction. As we could not find them in the literature in the exact form required for our applications, we include proofs. Results similar to the ones explained here are described in \cite{LM}*{\S~5.1, Theorems~6.7 and~6.11, Corollary~6.19} and \cite{BLV}*{Corollary~5.9}.

Let $\cT,\cV$ be tensor categories, where $\cV$ is symmetric with braiding $c$. Let $F\:\cT\to\cV$ be a tensor functor. In particular, such a functor is exact or equivalently faithful \cite{CEOP}*{Theorem~2.4.1}, and its extension to ind-completions $\tilde F\:\Ind(\cT)\to\Ind(\cV)$ has a right adjoint $\tilde G\:\Ind(\cV)\to\Ind(\cT)$ \cite{KS}*{Corollary~8.6.3}.

\begin{definition} A \emph{central section} of $F$ is a tensor functor $J\:\cV\to\cT$ together with isomorphisms $\gamma_{V,X}\: J(V)\o X\to X\o J(V)$ in $\cT$, natural and coherent with respect to tensor products in both $V\in\cV$ and $X\in\cT$ such that
$$
FJ \cong \Id_\cV 
\quad\text{ and }\quad
F(\gamma_{V,X}) = c_{V,FX} \quad
\text{for all }V,X.
$$
\end{definition}

Spelled out, the above coherences mean $\gamma_{V,\one}=\id$, $\gamma_{\one,X}=\id$, 
$$
\gamma_{V,X\o X'} = (\id\o\gamma_{V,X'})(\gamma_{V,X}\o\id) ,
\quad\text{and }\quad
\gamma_{V\o V',X} = (\gamma_{V,X}\o\id)(\id\o\gamma_{V',X}) .
$$

For instance, for a finite group $G$ and the forgetful functor $\Rep G\to\Vec$, a central section is given by the trivial-action functor, and similarly for categories of rational representations of algebraic groups. 

\begin{theorem} \label{thm:rec-1} There is a central section of $F$ if and only if there is a Hopf algebra $H$ in $\Ind(\cV)$ such that $\cT\simeq\Comod_\cV(H)$ as tensor categories with an equivalence that identifies $F$ with the forgetful functor.
\end{theorem}

\begin{proof} The ``if''-part follows from noting that a central section is given by the trivial-comodule functor that assigns any object $V\in\cV$ the trivial coaction of $H$ using its unit.

It remains to show the ``only if''-part. By Beck's (co)monadicity, $\Ind(\cT)\simeq\Comod_{\Ind(\cV)}(T)$ for the comonad $T:=\tilde F \tilde G$. By rigidity, we have projection formula isomorphisms in the sense of \cite{FLP}, see \cite{FLP}*{Corollary~3.20}. Hence, by \cite{BLV}*{Proposition~2.14}, applied in the dual situation, $T$ is a Hopf comonad. The central section equips $T$ with a monoidal coaugmentation compatible with the braiding in the sense of \cite{BLV}*{Corollary 5.9}. Hence, by that result, $T\simeq H\o -$, where $H:=T(\one)$ is a Hopf algebra in $\Ind(\cV)$. Passing to rigid (or compact) objects implies the assertion.
\end{proof}

While it is possible to upgrade \Cref{thm:rec-1} to include a description of the Hopf algebra $H$, we will derive a description that is more convenient in our set-up. For any object $X$ in a tensor category, let $T(X)=\bigoplus_{i\ge0} X^{\o i}$ be its tensor algebra in the ind-completion.

\begin{theorem} \label{thm:rec-2} Assume there is a central section of $F$. Assume also there is a self-dual object $W'\in\cT$ such that any object of $\cT$ is a subquotient of a finite direct sum of tensor powers of $W'$. Assume finally that there is a coalgebra structure $(\Delta,\eps)$ on $W:=F(W')$ and a morphism 
$$
\Delta'\:W'\to W'\o J(W)  \text{ in } \cT
\quad \text{ such that }
F(\Delta')=\Delta .
$$
Let $\psi$ be the image under $\tilde F$ of the adjunction mate $\psi'$ of the algebra map 
$$\tilde F (T(W'))\cong T(W)\to\one$$
given by $\eps$ on $W$, that is,
$$
\psi := \Big( T(W)\xrightarrow{\tilde F(\mathrm{unit})} \tilde F \tilde G (T(W)) \xrightarrow{\tilde F\tilde G((\eps^{\o i})_{i\ge0})} \tilde F\tilde G(\one) \Big)
\quad\text{ in } \Ind(\cV) .
$$
Then $\psi$ is an algebra morphism and the Hopf algebra $H$ as in \Cref{thm:rec-1} is isomorphic to $T(W)/\ker(\psi)$ as a bialgebra, with the algebra structure induced by that of the tensor algebra, and the coalgebra structure uniquely determined by the values of the structural morphisms on $W$ prescribed by $\Delta$ and $\eps$.
\end{theorem}

\begin{proof}
By \Cref{thm:rec-1}, we may assume $\cT=\Comod_\cV(H)$, $F$ is the forgetful functor, $G$ is the cofree comodule functor $X\mapsto {}^\bullet H\o X$, where ${}^\bullet H$ is the regular left $H$-comodule, $J$ is the trivial coaction functor, and the adjunction $(F,G)$ has unit
$$
\eta_X = \delta_X \: X\to {}^\bullet H\o F(X)=GF(X) ,
$$
where $\delta_X$ is the coaction of $H$ on $X$, for any $X\in\Comod_\cV(H)$. We set 
$$
\alpha := G(\eps)\eta_{W'}\: W'\to GF(\one)={}^\bullet H ,\quad
j := F(\alpha)\: W\to H .
$$
Let $\delta\:W\to H\o W$ be the coaction of $W'$ as an object in $\Comod_\cV(H)$. Then explicitly,
$$
j = (\id_H\o\eps) \delta .
$$

We claim that $j\:W\to H$ is a coalgebra morphism. Indeed, $\Delta'$ being a comodule morphism means
$$
 (\id_H\o F(\Delta')) \delta = (\delta\o\id_W) F(\Delta') .
$$
Using $F(\Delta')=\Delta$ and post-composing with $\id_H\o\eps\o\id_W$ yields
\begin{equation}
\delta = (j\o\id_W) \Delta  .
\label{eq:rec-delta}
\end{equation}
Now $\alpha$ being a comodule morphism and $j$ being its image in $\cV$, together with \Cref{eq:rec-delta} implies
$$
\Delta_H j = (\id_H\o j)\delta = (j\o j)\Delta ,
$$
where $\Delta_H$ is the comultiplication, i.e., the regular coaction of $H$. Post-composing \Cref{eq:rec-delta} with $\eps_H\o\eps$, we get
$$
\eps_H j = \eps  .
$$
We have established the claim.

The claim implies that there is a unique bialgebra morphism $\phi\:T(W)\to H$ given by $j$ on $W$, with the unique bialgebra structure on the algebra $T(W)$ that is given by $(\Delta,\eps)$ on $W$. Let $H'\subset H$ be the image of $\phi$, then \Cref{eq:rec-delta} shows that $\delta(W)\subset H'\o W$. As $W'$ generates $\cT$ in the sense described, this means $\delta_X(X)\subset H'\o X$ for all $X\in\cT$ and even all $X\in\Ind(\cT)$. In particular, for $X={}^\bullet H$,
$$
H = (\id_H\o \eps_H)\Delta_H(H)
 = (\id_H\o \eps_H)\delta_H({}^\bullet H)\subset H' ,
$$
i.e., $\phi$ is an epimorphism.

As $F$ is (strong) monoidal, the right adjoint $G$ has a canonical lax monoidal structure, and the unit and counit transformations are monoidal. In particular, $\psi'$ defined as the mate of an algebra morphism is an algebra morphism. Since it is given by $\alpha$ on $W'$, it is the unique such algebra morphism of the form $T(W')\to G(\one)$. In particular, $\psi=F(\psi')$ is the unique algebra morphism of the form $T(W)\to H$ that is $F(\alpha)=j$ on $W$. Hence, $\psi=\phi$, and $H\cong T(W)/\ker\psi$, as asserted.
\end{proof}

\section{Group-theoretical subcategories} \label{s:group-theoretical-subcats}

We start with the set-up of \Cref{bg:knop}:
Fix $\cA$ a small subobject-finite regular category. We assume that the terminal object $\one$ has no proper subobjects. Let $\delta$ be a degree function on $\cA$ valued in a fixed field $\kk$ of characteristic $0$, and let $\cT = \cT(\cA, \delta)$ be
the associated pseudo-tensor category. By restricting the collection of admissible relations, we obtain pseudo-tensor subcategories of $\cT$ which may admit an abelian envelope provided that $\cT$ admits one.

We fix a set $\cG$ of atoms of $\cA$, i.e., non-terminal objects admitting no proper product decomposition. We assume that the finite products in $\cA$ are chosen strictly associatively and unitally, and that its objects are exactly the chosen finite products
\[
    X_1\times\cdots\times X_n, \qquad n\geq 0,\quad X_i\in\cG,
\]
with the empty product equal to $\one$. Thus every object is literally such a product. The set $\cG$ and the products are fixed before choosing any subcategory; distinct members of $\cG$ may be isomorphic in $\cA$.

\begin{definition} \label{d:group-theoretical}
    A (not necessarily full) subcategory $\cC$ of $\Rel(\cA)$ containing all objects of $\cA$
    is called a \emph{group-theoretical subcategory} if $\cC$ is closed under

    \begin{enumerate}[label=(\arabic*)]
        \item forming evaluation and coevaluation relations: $\cC$ contains $X\times X \xleftarrow{(\id,\id)} X \xrightarrow{} \one$ and $\one \xleftarrow{} X \xrightarrow{(\id,\id)} X\times X$ for any $X \in \cG$.
        \item transpose: if $X \xleftarrow{f} r \xrightarrow{g} Y$ is a morphism in $\cC$, then so is its transpose $Y \xleftarrow{g} r \xrightarrow{f} X$, which we denote by $r^T$.
        \item ``tensor product'': if $X \xleftarrow{f} r \xrightarrow{g} Y$ and
            $X' \xleftarrow{f'} r' \xrightarrow{g'} Y'$ are morphisms in
            $\cC$, then so is $X \times X' \xleftarrow{f \times f'} r
            \times r' \xrightarrow{g \times g'} Y \times Y'$.
        \item subrelations: if $X \leftarrow r \to Y$ is a
            morphism in $\cC$ and $r' \hookrightarrow r$ is any subobject,
            then $X \leftarrow r' \rightarrow Y$ is also a morphism in
            $\cC$.
        \item automorphisms of tensor generators: For every $X \in \cG$ and all
            $f \in \mathrm{Aut}_{\cA}(X)$, the relation $X \xleftarrow{f} X \xrightarrow{=} X$ is contained in $\cC$.
    \end{enumerate}

    For a set $R$ of relations in $\Rel(\cA)$, we denote by $\gtcategory{R}$ the smallest group-theoretical subcategory of $\Rel(\cA)$ containing $R$ and call it the group-theoretical subcategory \emph{generated} by $R$. This smallest subcategory exists: $\Rel(\cA)$ is group-theoretical, and intersections of group-theoretical subcategories have the same fixed objects and satisfy conditions~(1)--(5). Thus $\gtcategory{R}$ is the intersection of all group-theoretical subcategories containing $R$.
\end{definition}

From now on, assume that $\cC$ is a group-theoretical subcategory
of $\Rel(\cA)$.
Using $\delta$, we can linearize $\cC$ to obtain a pseudo-tensor
category $\cC_\delta$ and a faithful linear monoidal
functor $\iota: \cC_\delta \to \cT$. Formally:

\begin{definition} \label{d:c-delta}
    Let $\cC_\delta^0$ be the subcategory of $\cT(\cA,\delta)$ on all objects of the form $[X]$, for $X\in\cA$, whose morphisms are given by linear combinations of relations in $\cC$. Let $\cC_\delta$ be the additive and idempotent completion of $\cC_\delta^0$.
\end{definition}

\begin{remark}
    For $X=X_1\times\cdots\times X_n$ with $X_i\in\cG$, the object $[X_n\times\cdots\times X_1]$ is a dual of $[X]$. Its evaluation and coevaluation are obtained by composing and tensoring those of the atoms. Thus $\cC_\delta$ is rigid. The pivotal structure inherited from $\cT$ makes $\cC_\delta$ spherical and we can talk about traces in $\cC_\delta$.
\end{remark}

\begin{remark}
    The assumption on products of atoms and conditions~(1)--(3) in \Cref{d:group-theoretical} make $\cC_\delta$ a pseudo-tensor subcategory of $\cT$. Conditions~(4) and~(5) ensure that $\cC_\delta$ ``sees enough'' of the structure of $\cA$ and $\cT$. In particular, condition~(4) implies that $\cC_\delta$ is semisimple whenever $\cT$ is (\Cref{t:c-delta-semisimple}), and is required in the construction of an abelian envelope of $\cC_\delta$ presented below.
\end{remark}

\begin{remark}
    Note that the categories $\Rel(\cA)$ and $\cT$ are invariant under replacing $\cA$ by an equivalent regular category, but this is not true for the class of group-theoretical subcategories of $\Rel(\cA)$. Indeed, if two objects in $\cA$ are isomorphic, they may be non-isomorphic in a subcategory of $\Rel(\cA)$. We explore this in \Cref{s:groups-in-rep-st}, where $\cA$ is taken to be equivalent to $\FinSet\op$, but not skeletal.
\end{remark}

\begin{theorem} \label{t:c-delta-semisimple}
    Assume that $\kk$ is algebraically closed. Then $\cC_\delta$ is semisimple for all the choices of $\delta$ that make $\cT(\cA,\delta)$ semisimple (as in \Cref{thm:knop}). These are exactly the semisimple choices of $\delta$ if $\cA$ is Barr-exact and Malcev.
\end{theorem}
\begin{proof}
    For an object $X \in \cA$, consider the set $\cC(X) \coloneqq \cC(\one, X) \subseteq \subobj(X)$ and a subobject $S \in \cC(X)$. By the assumptions on $\cA$, $\subobj(X)$ is a finite partial meet-semilattice with $S \wedge T = S\times_X T$ whenever the fiber product exists. Because $\cC$ is closed under coarsening, \Cref{l:moebius-function-restricts} asserts that $\subobj(S)$ and $\cC(X)$ are again partial meet-semilattices and the Möbius functions $\mu_{\subobj(S)}, \mu_{\cC(X)}$, and $\mu_{\subobj(X)}$ coincide on $\subobj(S)$, identified with the order ideal below $S$.

    Because $\one$ has no proper subobjects, we have $\End_{\cC_\delta}(\one) = \kk$. Denote by $\cN$ the tensor ideal of negligible morphisms in $\cC_\delta$. We have $\cN = 0$ if and only if for all $X \in \cA$, the bilinear form
    \[
        \beta_X: \cC_\delta(\one, [X]) \times \cC_\delta([X], \one) \to \End_{\cC_\delta}(\one) = \kk, \quad (f, g) \mapsto gf
    \]
    is non-degenerate for all $X \in \cA$. We can regard $\cC(X)$ as a basis of $\cC_\delta(\one, [X])$ and $\cC_\delta([X], \one)$, and in this basis, $\beta_X$ is given by $(S \hookrightarrow X, T \hookrightarrow X) \mapsto \delta(S \wedge T \tto \one)$. If the meet does not exist, we put $\delta(S \wedge T \tto \one) = 0$. Thus $\beta_X$ is non-degenerate if and only if
    \[
        \Omega_X \coloneqq \det(\delta(S \wedge T \tto \one)_{S,T \in \cC(X)}) \ne 0.
    \]
    Using \Cref{l:lindstroems-lemma} and the definition of the Möbius basis, we compute
    \[
        \Omega_X = \prod_{S \in \cC(X)} \; \sum_{\substack{T \in \cC(X)\\T \subseteq S}} \mu_{\cC(X)}(T, S) \delta(T \tto \one) = \prod_{S \in \cC(X)} \; \underbrace{\sum_{T \in \subobj(S)} \mu_{\subobj(S)}(T, S) \delta(T \tto \one)}_{\eqqcolon \omega_S}.
    \]
    Therefore, $\cN = 0$ if and only if for all $X \in \cA$ and all $S \in \cC(X)$ we have $\omega_S \ne 0$. Writing $X=\prod_i X_i$ with $X_i\in\cG$, the tensor product of their coevaluations gives a subobject isomorphic to $X$ in $\cC(\prod_i X_i^{\times 2})$. Thus this is equivalent to $\omega_X \ne 0$ for all $X \in \cA$.

    Now assume that $\cT(\cA, \delta)$ is semisimple. By \Cref{thm:knop}, $\omega_X \ne 0$ for all $X \in \cA$, so $\cN = 0$. Moreover, by \Cref{p:semisimplification} all nilpotent endomorphisms in $\cT(\cA, \delta)$ have trace 0. The traces in $\cC_\delta$ and $\cT(\cA, \delta)$ coincide, so \Cref{p:semisimplification} implies that $\cC_\delta$ is semisimple. Conversely, assume that $\cA$ is Barr-exact and Malcev and that $\cC_\delta$ is semisimple. By \Cref{p:semisimplification}, $\cN = 0$, so $\omega_X \ne 0$ for all $X \in \cA$. Thus $\cT(\cA, \delta)$ is semisimple by \Cref{thm:knop}.
\end{proof}

As explained in \Cref{thm:adjoint-ab-env}, a suitable monoidal adjunction can be used to conclude the existence of an abelian envelope of $\cC_\delta$.

In this section, we exhibit a condition under which the functor $\iota$ or its Ind-completion is part of such an adjunction; equivalently, under which the functor
\[
    \cT(\iota(-), \one): \cC_\delta\op \to \Vec \qquad\text{or}\qquad \tilde{\cT}(\tilde{\iota}(-), \one): \tilde{\cC}_\delta\op \to \Vec
\]
is representable, where the tilde denotes the Ind-completion.

\subsection{The rank assumption}

This subsection is independent of the choice of the degree function $\delta$ and we only work in the calculus of relations $\Rel(\cA)$. Recall that in $\cA$, pullbacks of surjective morphisms always exist, which explains the existence of the fiber products in this subsection.

\begin{notation}
    For a relation $X \xleftarrow{f} r \xrightarrow{g} Y$ in $\Rel(\cA)$,
    denote by $\topp(r)$ the subobject $\im(f) \hookrightarrow X$ and by
    $\bott(r)$ the subobject $\im(g) \hookrightarrow Y$. For any subobject
    $f: X \hookrightarrow Y$ in $\cA$, the relation $Y \xleftarrow{f} X
    \xrightarrow{f} Y$ is denoted by~$\beta_f$.
\end{notation}

Relations of the form $\beta_f$ are contained in any group-theoretical subcategory $\cC$ because they are subrelations of the identity relation on $Y$. In the following, we use juxtaposition to denote composition of relations (as for all other morphisms).

\begin{lemma} \label{l:topp-bott-composition}
    Let $r,s$ be relations in $\cC$ with $\bott(s)=\topp(r)$.

    \begin{enumerate}
        \item[(a)] We have $r \beta_{\topp(r)} = r = \beta_{\bott(r)} r$.
        \item[(b)] The composition $rs$ satisfies $\topp(rs) = \topp(s), \bott(rs) = \bott(r)$.
    \end{enumerate}
\end{lemma}
\begin{proof}
    Follows from the definition of composition.
\end{proof}

\begin{definition} \label{d:equivalence-of-tops}
    Consider the set $\cS = \bigsqcup_{X \in \cA} \subobj(X)$ of all
    subobjects of objects in $\cA$. We define an equivalence
    relation on $\cS$ by $f_1 \sim f_2$ iff there exists a relation $r \in
    \cC$ such that $\topp(r) = f_1, \bott(r) = f_2$.
\end{definition}

\begin{remark}
    The relation $\sim$ is reflexive because $\beta_f \in \cC(Y, Y)$ for
    all $f: X \hookrightarrow Y$. It is symmetric because relations in
    $\cC$ are closed under flipping the two legs. It is transitive by \Cref{l:topp-bott-composition}.
\end{remark}


\begin{assumption} \label{a:rank-object}
    We say that the group-theoretical subcategory $\cC \subseteq \Rel(\cA)$ \emph{has all rank objects} if for each subobject $f: X \hookrightarrow Y$, there exists an object $R_f \in \cC$ such that the following conditions hold:
    \begin{enumerate}
        \item[(A1)] there is a relation of the form $Y \xhookleftarrow{f}
            X \to R_f$ in $\cC$,
        \item[(A2)] for all relations of the form $Y \xleftarrow{f \eps}
             \cdot \to R_f$ in $\cC$ with surjective $\eps$, the
            morphism $\eps$ is an isomorphism.
    \end{enumerate}
\end{assumption}

\begin{remark}
    We can interpret the object $R_f$ as the ``rank'' of a subobject $f: X \hookrightarrow Y$ (with respect to the choice of $\cC$), i.e.~in some sense the smallest object for which there exists a relation $r \in \cC(Y, R_f)$ with $\topp(r) = f$. In the examples in \Cref{s:groups-in-rep-st,s:groups-in-rep-g-st,s:group-in-rep-glt}, this is how we verify \Cref{a:rank-object}.
\end{remark}

\begin{lemma} \label{l:equivalent-assumption}
    Assuming \Cref{a:rank-object}, the object $R_f$ can be chosen independently of the choice of representative $f$ in its equivalence class $[f]_{\sim}$.
\end{lemma}
\begin{proof}
    Consider a subobject $f': X' \hookrightarrow Y'$ with $f' \sim f$, i.e.~there exists a relation $s \in \cC(Y, Y')$ of the form $Y \xhookleftarrow{f} X \stackrel{}{\twoheadleftarrow} s \stackrel{}{\tto} X' \xhookrightarrow{f'} Y'$. We need to show that the conditions (A1), (A2) for $f$ imply the corresponding conditions (A1'), (A2') for $f'$.
    
    By (A1), there exists a relation of the form $Y \xhookleftarrow{f} X \to R_f$ in $\cC$. Composition with $s^T$ yields a relation of the form
    \[
        Y' \xhookleftarrow{f'} X' \twoheadleftarrow Z \to R_f.
    \]
    It will follow from (A2') that $Z \tto X'$ must be an isomorphism, which shows (A1').
    
    To prove (A2'), consider a relation $Y' \xhookleftarrow{f'} X' \stackrel{\eps}{\twoheadleftarrow} r \to R_f$. We need to show that $\eps$ is an isomorphism. Because we are working in a regular category, it suffices to show that $\eps$ is injective. To this end, let $p, q: A \rightrightarrows r$ be morphisms with $\eps p = \eps q$. We need to show that $p = q$.
    
    Consider the composition $rs$ and form the pullback of $s \tto X'$ and $\eps p = \eps q$. From the universal property of the fiber product $s\times_{X'} r$, we get two induced morphisms  $p', q'$:
    \[\begin{tikzcd}
        && s \times_{X'} A \arrow[bend right=60,swap]{lddd}{g} \arrow[two heads,bend left]{rrd}{h} \arrow[dashed,shift right,swap]{d}{p'} \arrow[dashed,shift left]{d}{q'} && \\[0.4cm]
        && s \times_{Y'} r = s \times_{X'} r \arrow{ldd} \arrow{rdd} \arrow[two heads]{d}{\zeta} && A \arrow[shift right,swap]{ldd}{p} \arrow[shift left]{ldd}{q} \\[0.4cm]
        && rs \arrow[hook]{d} && \\
        & s \arrow[two heads]{ld} \arrow[two heads]{rd} & Y \times R_f \arrow[dotted,bend left=15]{lldd}\arrow[dotted,bend right=15]{rrdd} & r \arrow[two heads,pos=.3]{ld}{\eps} \arrow{rdd}{m} & \\
        X \arrow[hook]{d}{f} && X' \arrow[hook]{d}{f'} && \\
        Y && Y' && R_f \\
    \end{tikzcd}\vspace{-.5cm}\]
    Write $\lambda: rs \hookrightarrow Y \times R_f \to Y$ and $\rho: rs \hookrightarrow Y \times R_f \to R_f$ for the left and right leg of the relation $r \circ s$. By (A2), $\lambda$ is injective, so 
    \[
        \lambda \zeta p' = f(s \tto X)g = \lambda \zeta q'
    \]
    implies $\zeta p' = \zeta q'$. Therefore
    \[
        mph = \rho\zeta p' = \rho\zeta q' = mqh.
    \]
    The morphism $h$ is a pullback of a surjection, hence epic, so $mp = mq$. In addition $f'\eps p = f'\eps q$, and since $f'\eps$ and $m$ are jointly injective, we get $p = q$.
\end{proof}

\begin{example} \label{e:terminal-occurs-as-rank}
    Consider the terminal object $\one$ of $\cA$ and $f = \id_\one \in \cS$. Then $R_f$ can be chosen to be $\one$. Indeed, condition~(A1) is satisfied because $\cC$ must contain $\one \xleftarrow{=} \one \xrightarrow{=} \one$. For condition~(A2), consider a relation $\one \xleftarrow{=} \one \stackrel{\eps}{\twoheadleftarrow} X \xrightarrow{g} \one$ in $\cC$. The map $X \to \one \times \one \cong \one$ induced by $\eps$ and $g$ is simply the unique morphism $\eps: X \to \one$ and it must be injective. Thus $\eps$ is an isomorphism.
\end{example}

\begin{proposition} \label{p:r-s-q}
    Assume \Cref{a:rank-object} and choose the rank objects as in \Cref{l:equivalent-assumption}.
    Let $f \sim f'$ be elements in $\cS$ and write $f': X
    \hookrightarrow Y$. Consider relations $r$ of the form $Y
    \xhookleftarrow{f'} X \xrightarrow{g} R_f$ and $s$ of the form $Y
    \xhookleftarrow{f'} X \xrightarrow{g'} R_f$ in $\cC$. Then there exists
    $R_f \hookleftarrow q \hookrightarrow R_f$ in $\cC$ such that $\bott(r)
    = \topp(q)$ and $s = qr$.
\end{proposition}
\begin{proof}
    Because $R_f = R_{f'}$, we may assume $f = f'$. Define $q \coloneqq s
    r^T$, i.e.~as the relation
    \[\begin{tikzcd}
        && X = X \times_Y X \arrow[equal]{ldd} \arrow[equal]{rdd} \arrow[two heads]{d}{\eta} && \\
        && q \arrow[hook]{d} && \\
        & X \arrow[swap]{ld}{g} \arrow[hook,swap,pos=.7]{rd}{f} & R_f \times R_f \arrow[dotted]{lld}\arrow[dotted]{rrd} & X \arrow[hook',pos=.7]{ld}{f} \arrow{rd}{g'} & \\
        R_f && Y && R_f. \\
    \end{tikzcd}\]
    By (A2), the two legs $q \hookrightarrow R_f \times R_f \to R_f$ are
    injective. This implies that $g$ and $g'$ equalize the
    same morphisms: We have $gp_1 = gp_2$ iff $\eta p_1 = \eta p_2$ iff
    $g'p_1 = g'p_2$ for all morphisms $p_1, p_2: A \to X$. The composition
    $qr = s r^T r$ is given by
    \[\begin{tikzcd}
        &&& X \times_{R_f} X \times_Y X \arrow[bend right=45]{llddd} \arrow[bend left=45]{rrddd} \arrow[bend left=60]{ddd} \arrow[two heads]{d} &&& \\
        &&& t \arrow[hook]{d} &&& \\
        &&& Y \times R_f \arrow[dotted]{llldd} \arrow[dotted]{rrrdd} &&& \\
        & X \arrow[swap,hook']{ld}{f} \arrow[swap]{rd}{g} && X \arrow{ld}{g} \arrow[hook,swap]{rd}{f} && X \arrow[hook']{ld}{f} \arrow{rd}{g'} & \\
        Y && R_f && Y && R_f \\
    \end{tikzcd}\]
    As for $q$, the pullback over the injection $f$ simplifies the diagram to
    \[\begin{tikzcd}
        && X \times_{R_f} X \arrow[swap]{ldd}{\pi_1} \arrow{rdd}{\pi_2} \arrow[two heads,pos=.4]{d}{\eps} && \\
        && t \arrow[hook]{d}{\iota} && \\
        & X \arrow[swap,hook']{ld}{f} \arrow[swap,pos=.7]{rd}{g} & Y \times R_f \arrow[dotted]{lld}\arrow[dotted]{rrd} & X \arrow[pos=.7]{ld}{g} \arrow{rd}{g'} & \\
        Y && R_f && R_f \\
    \end{tikzcd}
    \]
    The pullback $X \times_{R_f} X$ over $g$ is the kernel pair of $g$, so
    the two projections $\pi_1, \pi_2$ are epimorphisms. Thus the morphism $t \xhookrightarrow{\iota} Y \times R_f \to Y$ has image $f$. It is thus monic by (A2). Because image factorizations are unique, we get an
    isomorphism $\alpha: t \to X$ making the following diagram commute:
    \begin{equation} \label{e:def-alpha} \begin{tikzcd}[column sep=2cm,row sep=.4cm]
        & t \arrow[hook,start anchor={[yshift=-1pt]}]{rd} \arrow[dashed]{dd}{\alpha}[swap]{\cong} \arrow[hook]{r}{\iota} & Y \times R_f \arrow{d} \\
        X \times_{R_f} X \arrow[two heads]{ru}{\eps} \arrow[two heads]{rd}{\pi_1} && Y \\
        & X \arrow[hook]{ru}{f} & \\
    \end{tikzcd}\end{equation}
    To show that the relations $Y \hookleftarrow t \to R_f$ and $s$
    coincide (i.e.~are the same subobject of $Y \times R_f$), we show that the following diagrams commute:
    \[\begin{tikzcd}
        X \arrow[hook,swap]{d}{f} & t \arrow[swap]{l}{\alpha} \arrow[hook]{d}{\iota} \\
        Y & Y \times R_f \arrow{l}
    \end{tikzcd} \qquad\text{and}\qquad
    \begin{tikzcd}
        t \arrow[hook,swap]{d}{\iota} \arrow{r}{\alpha} & X \arrow{d}{g'} \\
        Y \times R_f \arrow[swap]{r}{\psi} & R_f
    \end{tikzcd}\]
    The first diagram commutes by the definition of
    $\alpha$~(see~\eqref{e:def-alpha}). For the second, note that we have
    $g\pi_2 = g\pi_1 =g\alpha\eps$. Because $g$ and $g'$ equalize the same morphisms, we get $g'\alpha\eps = g'\pi_2 = \psi\iota\eps$. Since
    $\eps$ is surjective, $g'\alpha = \psi\iota$ as desired.
\end{proof}

\subsection{The Möbius basis}

We now return to the original category $\cA$, without the formally adjoined initial object.

Subobjects of a fixed object are partially ordered by inclusion; we write
$f \geqslant f'$ if an injection $f'$ factors through another injection $f$. In particular, for $X, Y \in \cC$, the set $\cC(X, Y)$ of relations between $X$ and $Y$ is
partially ordered under the subobject order. In fact, $\cC(X, Y)$ is a
partial meet-semilattice. Since $\cA$ is subobject-finite, $\cC(X, Y)$ is finite and we can
consider the Möbius basis $\{x(r) \mid r \in \cC(X, Y)\}$ of $\cC_\delta([X],
[Y])$.

\begin{proposition} \label{p:moebius-property}
    Consider relations $X \hookleftarrow r \hookrightarrow Y$ in $\cC$ and $s \in
    \cC(W, X)$. Then $x(r) s = 0$ (in $\cC_\delta$) unless $\bott(s) \geqslant \topp(r)$. If
    $s$ is of the form $W \hookleftarrow s \hookrightarrow X$, then $x(r)
    x(s) = x(r) s = x(rs)$ if $\bott(s) = \topp(r)$, and $x(r) x(s) = 0$
    otherwise.
\end{proposition}
\begin{proof}
    Consider the partial meet-semilattice $L$ of subobjects of $X$ and its Möbius
    algebra $A(L)$. The general absorption property of the Möbius basis $\{x(f) \mid f \in L\} \subseteq A(L)$ states that
    \[
        x(f) g = \begin{cases}
            x(f) & \text{if } g \geqslant f, \\
            0 & \text{otherwise},
        \end{cases} \qquad\text{and}\qquad
        x(f) x(g) = \begin{cases}
            x(f) & \text{if } g = f, \\
            0 & \text{otherwise}.
        \end{cases}
    \]
    for all $f, g \in L$. These properties are preserved under the
    embedding of algebras $A(L) \hookrightarrow \cC_\delta([X], [X]), f \mapsto
    \beta_f$. Also, the lattices of coarsenings of $r$ and
    $\beta_{\topp(r)}$ are isomorphic, so $x(r) = r x(\beta_{\topp(r)})$.
    Thus $\bott(s) \not\geqslant \topp(r)$ implies $x(r) s = r
    x(\beta_{\topp(r)}) \beta_{\bott(s)} s = 0$.

    If $s$ is of the form $W \hookleftarrow s \hookrightarrow X$ with
    $\bott(s) = \topp(r)$, then we can identify $s \times_X r \cong r \cong
    s$, and the composition $rs$ is given by
    \[\begin{tikzcd}
        && s \times_X r \arrow[swap,]{ldd}{\cong} \arrow{rdd}{\cong} \arrow{d}{\cong} && \\
        && rs \arrow[hook]{d} && \\
        & s \arrow[hook']{ld} \arrow[hook]{rd} & W \times Y \arrow[dotted]{lld}\arrow[dotted]{rrd} & r \arrow[hook']{ld} \arrow[hook]{rd} & \\
        W && X && Y \\
    \end{tikzcd}\]
    Thus the meet-semilattices of coarsenings of $\beta_{\topp(r)}$, $r$,
    and $rs$ are all isomorphic, and $x(\beta_{\topp(r)})$ is an
    idempotent. We obtain
    \[
        x(r) x(s) = r x(\beta_{\topp(r)}) x(\beta_{\bott(s)}) s = r x(\beta_{\topp(r)}) s = x(r) s = x(rs).
    \]
    On the other hand, if $\bott(s) \neq \topp(r)$, then $x(r) x(s) = r
    x(\beta_{\topp(r)}) x(\beta_{\bott(s)}) s = 0$.
\end{proof}

\subsection{The representing object}

The object representing $\cT(\iota(-), \one)$ will have a direct summand for each equivalence class $[f]_{\sim} \in \cS/{\sim}$. In the following, we assume \Cref{a:rank-object} and choose one rank object $R_f$ for each $[f]_{\sim}$, which is possible by \Cref{l:equivalent-assumption}. We choose $R_{\id_\one}=\one$ as in \Cref{e:terminal-occurs-as-rank}.

For now, fix an object $X \in \cC$ and write $\cS_X = \{f \in \subobj(X)
\mid R_f = X\}$. The equivalence relation $\sim$ on $\cS$ restricts to an
equivalence relation $\approx$ on $\cS_X$.

\begin{lemma} \label{l:approx-equal}
    For $f \approx f' \in \cS_X$ with $f \geqslant f'$, we have $f = f'$.
\end{lemma}
\begin{proof}
    Write $f: A \hookrightarrow X$ and $f': A' \hookrightarrow X$. The
    condition $f \geqslant f'$ means that $f'$ factors as $A'
    \hookrightarrow A \xhookrightarrow{f} X$. Since $f \approx f'$, there
    exists a relation $X \leftarrow r \rightarrow X$ in $\cC$ with
    $\topp(r) = f$ and $\bott(r) = f'$. By (A2), the two morphisms $r \to
    X$ are actually injective, so $A \cong r \cong A'$. Because $\cA$ is
    subobject-finite, the injective endomorphism $A \cong A' \hookrightarrow A$ is
    an isomorphism and $f, f'$ coincide as subobjects.
\end{proof}

\begin{lemma}
    For $f \in \cS_X$, the number
    \[
        a_f \coloneqq \abs{\{r \in \cC(X, X) \mid \topp(r) = f\}}
    \]
    is independent of the choice of the representative $f$ in its
    equivalence class $[f]_{\approx}$.
\end{lemma}
\begin{proof}
    Let $f' \approx f$ and choose $r \in \cC(X, X)$ with $\topp(r) = f$ and
    $\bott(r) = f'$. By~(A2), we have $r^T r = \beta_f$ and $r r^T =
    \beta_{f'}$. Thus the map
    \[
        \{s \in \cC(X, X) \mid \topp(s) = f\} \longrightarrow
        \{s' \in \cC(X, X) \mid \topp(s') = f'\},
        \qquad s \longmapsto s r^T
    \]
    has an inverse given by $s \mapsto sr$.
\end{proof}

\begin{definition} 
    For $f \in \cS_X$, define
    \[
        e_f \coloneqq \frac{1}{a_f} \sum_{\substack{r \in \cC(X, X) \\
        \topp(r) \approx f}} x(r) \in \cC_\delta([X], [X]),
    \]
    where the $x(r)$ are taken in the Möbius basis of $\cC_\delta([X], [X])$.
\end{definition}

\begin{lemma}
    For $f \in \cS_X$, the element $e_f$ is an idempotent in $\cC_\delta([X],
    [X])$.
\end{lemma}
\begin{proof}
    By~(A2), all $r \in \cC(X, X)$ with $\topp(r) \approx f$
    are of the form $X \hookleftarrow r \hookrightarrow X$. To see that
    $e_f$ is an idempotent, we compute
    \[\begin{aligned}
        e_f^2 &= \frac{1}{a_f^2} \sum_{\substack{r \in \cC(X, X) \\
        \topp(r) \approx f}} \; \sum_{\substack{s \in \cC(X, X) \\ \topp(s)
        \approx f}} x(r) x(s) \\
        &\stackrel{(*)}{=} \frac{1}{a_f^2} \sum_{\substack{r \in \cC(X, X)
        \\ \topp(r) \approx f}} \; \sum_{\substack{s \in \cC(X, X) \\
        \bott(s) = \topp(r)}} x(rs) \\
        &= \frac{1}{a_f^2} \sum_{\substack{t \in \cC(X, X) \\ \topp(t)
        \approx f}} \; \sum_{\substack{rs = t \\ \bott(s) = \topp(r)}} x(t)
        \\
        &\stackrel{(**)}{=} \frac{1}{a_f^2} \; \sum_{\substack{t \in \cC(X,
        X) \\ \topp(t) \approx f}} \; a_f x(t) \\
        &= e_f,
    \end{aligned}\]
    where the equality $(*)$ follows from~\Cref{p:moebius-property} and
    $(**)$ follows from the fact that for fixed $t \in \cC(X, X)$ with
    $\topp(t) \approx f$, we have
    \[\begin{aligned}[b]
        &\abs{\{(r, s) \mid r, s \in \cC(X, X), rs = t, \topp(r) =
        \bott(s)\}} \\
        =\;&\abs{\{(r, r^Tt) \mid r \in \cC(X, X), \bott(r) = \bott(t)\}}
        = a_f.
    \end{aligned} \qedhere\]
\end{proof}

\begin{lemma} \label{l:ef-properties}
    Let $f \in \cS_X$ and $W \in \cC$. The map
    \[
        \cC_\delta([W], \im(e_f)) \to \cT([W], \one)
    \]
    given by postcomposition with the relation $X \xleftarrow{=} X \to \one$ is injective. Its image has a basis of the form
    \[
        \{f' + \text{lower order terms} \mid f' \in \subobj(W), f' \sim f\},
    \]
    where ``lower order'' refers to the subobject order on $\subobj(W)$ and we identify
    $\cT([W], \one) \cong \sspan_{\kk} \subobj(W)$.
\end{lemma}
\begin{proof}
    Consider $r \in \cC(X, X)$ and $s \in \cC(W, X)$ with $f \approx
    \bott(s) \approx \topp(r)$. By~\Cref{l:approx-equal}, the condition
    $\bott(s) \geqslant \topp(r)$ implies $\bott(s) = \topp(r)$. Thus $x(r)
    s = 0$ unless $\bott(s) = \topp(r)$ by~\Cref{p:moebius-property}, and
    we get
    \begin{equation}\label{e:efs}
        e_f s = \frac{1}{a_f} \sum_{\substack{r \in \cC(X, X) \\ \topp(r)
        \approx f}} x(r) s = \frac{1}{a_f} \sum_{\substack{r \in \cC(X, X)
        \\ \bott(s) = \topp(r)}} x(r) s.
    \end{equation}
    If $q \in \cC(X, X)$ satisfies $\topp(q) \approx f$ as subobjects of $X$, then again
    by~\Cref{p:moebius-property},
    \begin{equation}\label{e:efq}
        e_fq = \frac{1}{a_f} \sum_{\substack{r \in \cC(X, X) \\ \topp(r)
        = \bott(q)}} x(r) q = \frac{1}{a_f} \sum_{\substack{r \in \cC(X, X)
        \\ \topp(r) = \bott(q)}} x(rq) = e_f\beta_{\topp(q)}.
    \end{equation}
    For general $s \in \cC(W, X)$, we have
    \[\begin{aligned}
        e_f s
        &= \frac{1}{a_f} \sum_{\substack{r \in \cC(X, X) \\ \topp(r)
        \approx f}} x(r) s = \frac{1}{a_f} \sum_{\substack{r \in \cC(X, X)
        \\ f \approx \topp(r) \leqslant \bott(s)}} x(r) s
        = \frac{1}{a_f} \sum_{\bott(s) \geqslant \hat{f} \approx f} \;
        \sum_{\substack{r \in \cC(X, X) \\ \topp(r) = \hat{f}}} x(r) s \\
        &= \sum_{\bott(s) \geqslant \hat{f} \approx f} \; \frac{1}{a_f}
        \sum_{\substack{r \in \cC(X, X) \\ \topp(r) = \hat{f}}} x(r)
        (\beta_{\hat{f}}s)
        \stackrel{\eqref{e:efs}}{=} \sum_{\bott(s) \geqslant \hat{f}
        \approx f} \; e_{f} (\beta_{\hat{f}}s)
    \end{aligned}\]
    Therefore $\cC_\delta([W], \im(e_f))$ is spanned by $\{e_fs \mid s \in
    \cC(W, X), \bott(s) \approx f\}$. For each subobject $f \sim f' \in
    \subobj(W)$, \Cref{l:equivalent-assumption} implies that there exists a relation $s_{f'} \in \cC(W, X)$ with $\topp(s_{f'}) =
    f'$. If another $s \in \cC(W, X)$ satisfies $\topp(s) = f'$, then
    by~\Cref{p:r-s-q} there exists $q \in \cC(X, X)$ with $s = qs_{f'}$
    and $\topp(q) = \bott(s_{f'})$. Thus $e_f s = e_f q s_{f'} = e_f
    \beta_{\bott(s_{f'})}s_{f'} = e_f s_{f'}$ by~\eqref{e:efq}.

    In total $\{e_f s_{f'} \mid f' \in \subobj(W), f \sim f'\}$ is a spanning
    set. The map
    \[
        \cC_\delta([W], \im(e_f)) \to \cT([W], \one)
    \]
    sends $e_f s_{f'}$ to $f' + \text{lower order terms}$, and elements of this form are
    linearly independent.
\end{proof}

Recall that we denote the $\Ind$-completion of categories and functors by a tilde.

\begin{definition}
    We define the $\Ind$-object
    \[
        H_\cC \coloneqq \bigoplus_{X \in \cC} \; \bigoplus_{[f]_{\approx} \in
        \cS_X/{\approx}} \im(e_f) \in \tilde{\cC}_\delta.
    \]
\end{definition}

\begin{proposition} \label{p:representing-object}
    Consider the morphism $p: \tilde{\iota}(H_\cC) \to \one$ given on $\im(e_f)$ for $X \in \cC$ and $[f]_{\approx} \in \cS_X/{\approx}$ by the relation $X \xleftarrow{=} X \to \one$. Then the linear map
    \[
        \tilde{\cC}_\delta([W], H_\cC) \longrightarrow \tilde{\cT}([W], \one),
        \qquad s \longmapsto p \tilde{\iota}(s)
    \]
    is an isomorphism for all $W \in \cC$.
\end{proposition}
\begin{proof}
    The map is a direct sum of the maps $\cC_\delta([W], \im(e_f)) \to
    \cT([W], \one)$ from~\Cref{l:ef-properties} for $X \in \cC$
    and $[f] \in \cS_X/{\approx}$. These maps are injective and $\cT([W], \one)$ is precisely the direct sum of their images.
\end{proof}

\begin{corollary} \label{c:representing-object}
    The functor $\tilde{\cT}(\tilde{\iota}(-), \one): \tilde{\cC}_\delta\op \to \mathrm{Vec}$ is represented by the object $H_\cC$.
\end{corollary}
\begin{proof}
    The isomorphisms from \Cref{p:representing-object} are natural because they are given by post-composition with a fixed morphism. They extend to an isomorphism of functors on the respective Karoubi completions and $\Ind$-completions.
\end{proof}

\begin{lemma} \label{l:direct-summand-one}
    The object $H_\cC$ has a direct summand isomorphic to the unit object $\one \in \cC_\delta$.
\end{lemma}
\begin{proof}
    By~\Cref{e:terminal-occurs-as-rank}, the terminal object $\one \in \cA$ occurs as a rank object $R_{\id_\one}$. We have $e_{\id_\one} = \id_\one$ and $\one = \im(e_{\id_\one}) \in \cC_\delta$ is a direct summand in $H_\cC$.
\end{proof}

\begin{proposition} \label{prop:adjoint} For any group-theoretical category $\cC$ in $\Rel(\cA)$ satisfying \Cref{a:rank-object}, the inclusion functor $\iota$ has an ind-right adjoint such that the adjunction has a split unit.
\end{proposition}

\begin{proof}
    The functor $\iota: \cC_\delta \to \cT$ is dominant, because $\cC_\delta$ contains all objects $[X]$ with $X\in\cA$ and every object of $\cT$ is a direct summand of a finite direct sum of such objects. Therefore, by \Cref{l:ind-adjoint} and \Cref{c:representing-object}, the functor $\Ind(\iota): \Ind(\cC_\delta) \to \Ind(\cT)$ has a right adjoint $G$ sending $\one \mapsto H_\cC$. By \Cref{l:direct-summand-one} and the argument in the proof of \Cref{thm:adjoint-ab-env}, the adjunction unit splits naturally on $\cC_\delta$. Since $\Ind(\iota)$ preserves compact objects, $G$ preserves filtered colimits, so this splitting extends to the $\Ind$-completion.
\end{proof}

\begin{theorem} \label{t:group-theoretical-ab-env}
    Let $\cC$ be a group-theoretical category in $\Rel(\cA)$ such that the inclusion functor $\iota$ has an ind-right adjoint $G$ with a split unit (e.g., $\cC$ satisfies \Cref{a:rank-object}). Assume that $\cT$ has a global splitting object of the form $\iota(S)$ for some $S \in \cC_\delta$. Then $S$ is a global splitting object in $\cC_\delta$, $\cT$ and $\cC_\delta$ have abelian envelopes with enough projectives, and $\iota$ induces an (exact) tensor functor $\iota\ab\:\cC_\delta\ab\to\cT\ab$ between them.
\end{theorem}
\begin{proof}
    This is \Cref{thm:adjoint-ab-env}: the split unit makes $\one$ a direct summand of $G(\one)$. Under \Cref{a:rank-object}, the required adjunction is supplied by \Cref{prop:adjoint}.
\end{proof}

Now we can apply Tannakian reconstruction as laid out in \Cref{ss:tannaka-duality} and identify $\cC_\delta\ab$ with the category of comodules over the image of $G(\one)$ internal to $\cT\ab$:

\begin{theorem} \label{t:group-theoretical-comod}
    In the situation of \Cref{t:group-theoretical-ab-env}, the tensor functor $\iota\ab: \cC_\delta\ab \to \cT\ab$ factors as
    \[
        \cC_\delta\ab \xrightarrow{\cong} \Comod_{\cT\ab}(H) \xrightarrow{\text{forget}} \cT\ab,
    \]
    where $H=\Ind(\iota\ab)(G(\one))\in\Ind(\cT\ab)$, with the Hopf algebra structure reconstructed from the monoidal adjunction and a central section. Here we view $G(\one)$ as an object of $\Ind(\cC_\delta\ab)$ via the embedding of $\cC_\delta$. Under \Cref{a:rank-object}, it is the explicit object $H_\cC$ from \Cref{p:representing-object}.
\end{theorem}
\begin{proof}
    By \Cref{thm:rec-1}, it suffices to construct a central section of $\iota\ab: \cC_\delta\ab \to \cT\ab$. For any $X \in \cG$, define $\tilde X:=X^{\times 2}$ and consider the relation
    \[
        \alpha_X \coloneqq \left(\tilde X \xleftarrow{\Delta_X} X \xrightarrow{\id} X\right) \in \Rel(\cA)(\tilde X, X),
    \]
    where $\Delta_X$ is the diagonal map.
    For a general object $X \in \cA$, choose an expression $X=\prod_i X_i$ with $X_i\in\cG$. Set $\tilde X:=\prod_i \tilde X_i$, and let $\alpha_X\in\Rel(\cA)(\tilde X, X)$ be the product of all $\alpha_{X_i}$.
    For a relation $X \xleftarrow{f} r \xrightarrow{g} Y$ in $\Rel(\cA)$, the relation $J(r) \coloneqq \alpha_Y^T r \alpha_X$ is a coarsening of
    \[
        \tilde X \xleftarrow{(\prod_i \Delta_{X_i})\pi_X} X \times Y \xrightarrow{(\prod_j \Delta_{Y_j})\pi_Y} \tilde Y,
    \]
    using decompositions $X=\prod_i X_i$ and $Y=\prod_j Y_j$; therefore, $J(r)$ is contained in $\cC$. Note that in $\cT$, $\alpha_X\alpha_X^T = \id_{[X]}$ and $e_X \coloneqq \alpha_X^T\alpha_X$ is an idempotent. Put $J([X]) \coloneqq \im(e_X)$. The tensor constraint $J([X])\o J([Y])\to J([X\times Y])$ is induced by $\alpha_{X\times Y}^T(\alpha_X\o\alpha_Y)$. This relation and its transpose belong to $\cC$ by the same coarsening argument and induce inverse isomorphisms. Coherence follows after applying the faithful functor $\iota$. These assignments extend to a linear monoidal functor $J: \cT \to \cC_\delta$ and the $\alpha_X$ assemble to an isomorphism of monoidal functors $\alpha: \iota J \cong \id_\cT$. By the universal property of abelian envelopes, $J$ extends to a tensor functor $J\ab: \cT\ab \to \cC_\delta\ab$ with $\iota\ab J\ab \cong \id_{\cT\ab}$.

    For objects $X, Y \in \cA$, the relation
    \[
        \gamma_{X,Y} \coloneqq \left(\tilde X \times Y \longleftarrow X \times Y \longrightarrow Y \times \tilde X\right) \in \cT(\im(e_X) \o [Y], [Y] \o \im(e_X))
    \]
    is a coarsening of $\bigl(\bigotimes_i\ev_{[X_i]}\bigr) \o \id_{[Y]} \o \bigl(\bigotimes_i\coev_{[X_i]}\bigr)$, using the chosen decomposition $X=\prod_i X_i$, and therefore contained in $\cC_\delta$. It satisfies
    \[
        (\id_{[Y]} \o \;\alpha_X) \iota(\gamma_{X,Y}) (\alpha_X^T \otimes \id_{[Y]}) = c_{[X],[Y]}
    \]
    where $c$ is the symmetric braiding on $\cT$ given by swapping factors. Because the tensor product in $\cC_\delta\ab$ is exact, it commutes with direct sums and cokernels. Thus the definition of $\gamma$ extends to all $X \in \cT\ab, Y \in \cC_\delta\ab$ and makes $J\ab$ into a central section of $\iota\ab$.

    To identify the Hopf algebra, put $K=G(\one)$ and view it as an object of $\Ind(\cC_\delta\ab)$. By \cite{flake2026monoidaladjunctionsabelianenvelopes}*{Theorem~6.1}, every object $A \in \cC_\delta\ab$ has a presentation $P_1\to P_0\to A\to0$ by objects $P_i\in\cC_\delta$. Applying the adjunction isomorphisms $\tilde{\cC}_\delta(P_i,K)\cong\cT(\iota(P_i),\one)$ and taking kernels yields
    \[
        \Ind(\cC_\delta\ab)(A,K) \cong \cT\ab(\iota\ab(A),\one).
    \]
    Thus $K$ is the value at $\one$ of the right adjoint of $\Ind(\iota\ab)$, and \Cref{thm:rec-1} gives the asserted identification of $H$. Under \Cref{a:rank-object}, \Cref{c:representing-object} identifies $K$ with $H_\cC$.
\end{proof}

\subsection{Relation to oligomorphic groups} \label{ss:oligo}

We briefly recall from \cite{knoptooligo} how Knop categories fit into the more general construction of pseudo-tensor categories as ``permutation modules'' of oligomorphic groups \cite{oligotensor}. In this subsection, assume in addition that $\cA$ is Barr-exact.

\newcommand\Regg{\textbf{(Reg+)}\xspace}
For the subgroup construction below, we write \Regg for the following conditions: regular, Barr-exact, Malcev, subobject-finite, and with terminal object having no proper subobjects.

\begin{construction}[see \cite{knoptooligo}*{\S~4}] \label{c:g}
    The category $\Pro(\cA)$ of pro-objects over the regular category $\cA$ admits an object $\Omega$ such that for all $X \in \cA$, the set $\fB(X) \coloneqq \Surj_{\Pro(\cA)}(\Omega, X)$ is non-empty and a transitive $\Aut(\Omega)$-set. Such an object $\Omega$ is called a \emph{universal homogeneous pro-object} of $\cA$. The topology defined by the $\fB(X)$ makes $G \coloneqq \Aut(\Omega)$ into a pro-oligomorphic group.
\end{construction}

\begin{proposition} \label{p:bx-finite-cover}
    Every smooth transitive $G$-set has a finite cover by a $G$-set of the form $\fB(X)$ for some $X \in \cA$.
\end{proposition}
\begin{proof}
    Follows directly from \cite{knoptooligo}*{Proposition~4.4}.
\end{proof}

\begin{construction}[see \cite{knoptooligo}*{\S~6}] \label{c:mu}
    The degree function $\delta$ gives rise to a measure $\mu$ on $G$ in the following way: For a surjection $e: X \to Y$ in $\cA$, define
    \begin{equation} \label{e:define-mu}
        \mu(\fB(X) \to \fB(Y)) \coloneqq \sum_{\substack{Z \in \subobj(X),\\e(Z)=Y}} \mu_{\subobj(X)}(Z, X) \delta(e|_Z).
    \end{equation}
    Note that this is precisely Knop's element $\omega_e$ (see \cite{knop}*{\S~8}). Since we are in characteristic~0, \Cref{p:bx-finite-cover} implies that the assignments~\eqref{e:define-mu} extend uniquely to a measure on all finitary smooth $G$-sets.
\end{construction}

In the following, let $(G, \mu)$ be as in \Cref{c:g,c:mu}. Part~(a) of the following result was pointed out to us by Thorsten Heidersdorf and Andrew Snowden.

\begin{proposition} \label{p:knoptooligo}
    (a) There is an equivalence of pseudo-tensor categories
    \[
        \Phi: \cT(\cA, \delta) \longrightarrow \uPerm(G, \mu)\kar
    \]
    sending any object $[X]$ for $X \in \cA$ to $\fA(X) \coloneqq \bigsqcup_{Y \in \subobj(X)} \fB(Y)$ and any relation $r \hookrightarrow X \times Y$ in $\cA$ to the indicator function $1_{\fA(r)}$, regarded as a $\fA(X) \times \fA(Y)$ matrix.

    (b) If $\cA$ is additionally Malcev, then in the equivalent pseudo-tensor categories from part~(a), nilpotent endomorphisms have trace zero.
\end{proposition}
\begin{proof}
    (a) By~\cite{knoptooligo}*{Theorem~7.5}, the functor $\Phi$ exists as a fully faithful monoidal functor whose image is generated by $G$-sets of the form $\fB(X)$. For any smooth transitive $G$-set $T$, \Cref{p:bx-finite-cover} yields a finite cover $\pi: \fB(X) \tto T$ of degree $d$. The endomorphism $d^{-1}\pi^*\pi_*$ is an idempotent whose image is isomorphic to $T$ (recall that we are in characteristic 0). Thus $\Phi$ is essentially surjective.

    (b) By \cite{knop}*{Theorem~6.1}, $\cT(\cA,\delta)$ admits a linear symmetric monoidal functor to a semisimple tensor category. As the trace is preserved by such a functor, and as the trace of a nilpotent endomorphism is zero in any tensor category, this implies the assertion.
\end{proof}

The above equivalence is particularly useful when combined with the following description of categories corresponding to open subgroups of $G$. For $X\in\cA$, let $\cA//X$ be the category of surjections to $X$ in $\cA$, with commutative triangles as morphisms.

\begin{proposition} \label{p:subgroup-of-G}
    Assume that $\cA$ satisfies \Regg, and let $(G,\mu)$ be associated to $(\cA,\delta)$ as above.
    \begin{enumerate}
        \item[(a)] For every $X\in\cA$, the category $\cA//X$ satisfies \Regg. The rule $\delta_X(e)\coloneqq\delta(e)$ defines a degree function on $\cA//X$.
        \item[(b)] For every open subgroup $U$ of $G$, there is a finite-index open subgroup $V$ of $U$ and an object $X\in\cA$ such that $\fB(X)\cong G/V$ as $G$-sets and the pair $(V,\mu|_V)$ corresponds to the pair $(\cA//X,\delta_X)$ via \Cref{c:g,c:mu}.
    \end{enumerate}
\end{proposition}

\begin{proof}
    (a) Regularity of $\cA//X$ is explained in \cite{knop}*{\S~2, Example~4}, with image factorizations and pullbacks of surjections computed in $\cA$. A subobject in $\cA//X$ of a surjection $q\:Y\tto X$ is precisely a subobject $Z\subset Y$ in $\cA$ such that $q|_Z$ is surjective. Thus $\cA//X$ is subobject-finite, and its terminal object $\id_X$ has no proper subobjects. A reflexive relation in $\cA//X$ is reflexive in $\cA$, hence an equivalence relation. Its quotient in $\cA$ inherits a surjection to $X$. This proves exactness and the Malcev property for $\cA//X$ (see also \cite{knop}*{Appendix~A, Example~10}). The axioms of a degree-function for $\delta_X$ follow from those for $\delta$.

    (b) By \Cref{p:bx-finite-cover}, we can choose a finite cover $\pi:\fB(X)\tto G/U$. Let $\Omega$ be a universal homogeneous pro-object of $\cA$ and pick an element $p:\Omega\tto X$ of $\fB(X)$ with $\pi(p)=U$, the identity coset in $G/U$. Its stabilizer subgroup $V\subseteq G$ is open and contained in $U$. As $\fB(X)$ is a transitive $G$-set, $\fB(X)\cong G/V$, and the preimage of $U$ under $\pi$ is isomorphic to $U/V$, so $[U:V]<\infty$.

    Set $\cD:=\cA//X$. The quotients of $\Omega$ through which $p$ factors form a cofinal system. Let $\Omega'$ be the pro-object of $\cA//X=\cD$ defined by this system. We claim that $\Omega'$ is a universal homogeneous pro-object of $\cD$ in the sense of \Cref{c:g} with $\Aut(\Omega')\cong V$ as topological groups.
    
    Indeed, for any $q:Y\tto X$ in $\cD$, 
    \[
        \fB_\cD(q):=\Surj_{\Pro(\cD)}(\Omega',q)=\{f\in\fB(Y)\mid q f=p\}.
    \]
    As $\fB(Y)$ and $\fB(X)$ are non-empty transitive $G$-sets, the map of $G$-sets $f\mapsto qf$ between them is surjective, which implies that $\fB_\cD(q)$ is non-empty. For $f_1,f_2\in\fB_\cD(q)$ take $g\in G$ such that $f_2=f_1g$, then $p=q f_2=q f_1 g=pg$, i.e., $g\in V$ and $\fB_\cD(q)$ is transitive as a $V$-set. A similar argument shows that $V$ identifies with $\Aut(\Omega')$ as topological groups. In particular, to compare the neighborhood bases, note that for any $f:\Omega\tto Y$, the quotient of $\Omega$ given by the coimage of $(p,f):\Omega\to X\times Y$ has stabilizer $V\cap\operatorname{Stab}_G(f)$ and defines a point of $\fB_\cD(q)$ for some $q\in\cD$. This implies that we can identify basic neighborhoods in the subspace topology with stabilizer subgroups of points in the sets $\fB_\cD(q)$.
    We have shown the claim.

    As observed in the proof of (a), subobjects of $q:Y\tto X$ in $\cD$ are precisely the subobjects of $Y$ in $\cA$ which surject onto $X$ under $q$, and the intervals above these subobjects agree in the two subobject posets. For a surjection $e:(q:Y\tto X)\tto(r:Z\tto X)$ in $\cD$, every $T\subseteq Y$ with $e(T)=Z$ already surjects onto $X$ under $q$. Thus the sums in \Cref{e:define-mu} agree for $e$ in $\cD$ and in $\cA$. Since $\fB_{\cD}(e)$ is obtained from $\fB(e)$ by taking the fibers over $p$, this identifies the measure associated to $(\cD,\delta_X)$ with $\mu|_V$.
\end{proof}

\begin{lemma} \label{l:cdelta-contains-bx}
    In the situation of \Cref{p:knoptooligo}(a), let $\cC$ be a group-theoretical subcategory in $\Rel(\cA)$. Then $\Phi(\im(x(\id_X)))\cong\fB(X)$ for all $X\in\cA$. In particular, $\Phi(\cC_\delta)\subset\uPerm(G,\mu)\kar$ contains all objects of the form $\fB(X)$ for $X \in \cA$.
\end{lemma}
\begin{proof}
    Let $X \in \cA$ and $U \in \subobj(X)$ be a subobject. Regarded as a relation $X \hookleftarrow U \hookrightarrow X$ (i.e.~a morphism in $\cC_\delta$), we have
    \[
        \Phi(U) = 1_{\fA(U)} = \sum_{V \subseteq U} 1_{\fB(V)},
    \]
    where $\fB(V) \subseteq \fA(X) \subseteq \fA(X) \times \fA(X)$ via the diagonal embedding. Now consider the Möbius basis $x(U)_{U \in \subobj(X)}$. Möbius inversion yields $\Phi(x(U)) = 1_{\fB(U)}$. Taking $U = X$, we obtain
    \[
        \Phi(\im(x(X))) \cong \im(\Phi(x(X))) \cong \im(1_{\fB(X)}) \cong \fB(X),
    \]
    as asserted. The last assertion follows as $\cC$ is closed under coarsenings, in particular, those used to define $x(X)=x(\id_X)$.
\end{proof}

The following result establishes one of the hypotheses of \Cref{t:group-theoretical-ab-env} using the concept of quasi-regular measures.

\begin{proposition} \label{p:quasi-regular-splitting}
    Let $\cA,\delta$ be as above, assume $\cA$ is additionally Malcev. Let $(G, \mu)$ be the associated pro-oligomorphic group with a measure, as in \Cref{c:g,c:mu}, and assume $\mu$ is quasi-regular. Then there is a global splitting object $S$ in $\cT(\cA,\delta)$ that is contained in any subcategory $\cC_\delta$ corresponding to a group-theoretical subcategory $\cC$ of $\Rel(\cA)$.
\end{proposition}
\begin{proof}
    If $\mu$ is quasi-regular, there exists an open subgroup $U \subseteq G$ such that the restriction of $\mu$ to $U$ is regular. By \Cref{p:subgroup-of-G}(b), we can choose an open subgroup $V\subseteq U$ of finite index and $X\in\cA$ such that $\fB(X)\cong G/V$ and $(V,\mu|_V)$ is the pair associated to $(\cA//X,\delta_X)$. The restriction of $\mu$ to $V$ is still regular.

    By \Cref{p:subgroup-of-G}(a), the category $\cA//X$ satisfies \Regg. Hence, by \Cref{p:knoptooligo}(a) and~(b), nilpotent endomorphisms in $\uPerm(V,\mu|_V)\kar$ have trace zero.

    Hence, by the proof of \cite{oligotensor}*{Theorem~13.13}, using also \cite{oligotensor}*{Remark~11.5} if $G$ is not first-countable, $G/V$ is a global splitting object in $\uPerm(G, \mu)\kar$. Note that this proof uses what is called Property (P) to invoke \cite{oligotensor}*{Theorem~13.2}, but it is explained right after this result that it suffices to assume that nilpotent endomorphisms have trace zero, which we have just established for $V$.

    Set $S\coloneqq\im(x(\id_X))$ in $\cT(\cA,\delta)$. By  \Cref{l:cdelta-contains-bx}, $\Phi(S)\cong\fB(X)\cong G/V$. The monoidal equivalence in \Cref{p:knoptooligo}(a) therefore shows that $S$ is a global splitting object in $\cT(\cA,\delta)$ and for every group-theoretical subcategory $\cC\subseteq\Rel(\cA)$, $S$ lies in $\cC_\delta$, again by \Cref{l:cdelta-contains-bx}.
\end{proof}

\subsection{Regular categories with one atom} \label{ss:one-atom}

Consider a subobject-finite regular category $\cA$ where every object is a direct product $X^n$ for some fixed non-terminal object $X$ and $n \in \NN_0$. Write $\one \coloneqq X^0$ for the terminal object. We describe a special group-theoretical subcategory $\cC \subseteq \Rel(\cA)$ which does not satisfy \Cref{a:rank-object} but still yields an abelian envelope.

For a degree function $\delta$ on $\cA$, denote by $t$ the degree of the surjective morphism in the image factorization of $X \to \one$. Consider the relation
\[
    r \coloneqq \left(X^2 \xleftarrow{=} X^2 \to \one\right).
\]
Its rotation $\rho\coloneqq(X\xleftarrow{\pi_1}X^2\xrightarrow{\pi_2}X)$ satisfies $\rho^2 = t\rho$, so $t^{-1}\rho$ is an idempotent whenever $t \ne 0$.

\begin{lemma} \label{l:r-is-special}
    Let $\cC \subseteq \Rel(\cA)$ be a group-theoretical subcategory.
    \begin{enumerate}
        \item[(a)] The group-theoretical subcategory $\gtcategory{r} \subseteq \Rel(\cA)$ contains precisely the relations $X^n \leftarrow X^b \rightarrow X^m$ with $n + m$ even.
        \item[(b)] If $\cC$ contains $r$, then $\cC = \gtcategory{r}$ or $\cC = \Rel(\cA)$.
    \end{enumerate}
\end{lemma}
\begin{proof}
    For (a), note that $\gtcategory{r}$ contains the relation $X^{2n} \xleftarrow{=} X^{2n} \to \one$ for all $n \in \NN_0$. By coarsening and rotation, we see that $\gtcategory{r}$ contains all relations $X^n \leftarrow X^b \rightarrow X^m$ with $n + m$ even. The converse follows from the fact that $\gtcategory{r}$ is generated by such relations under tensor product, composition, reflection and coarsening, all of which preserve the parity.

    For (b), assume that $\cC$ also contains a relation $X^n \leftarrow X^b \rightarrow X^m$ with $n + m$ odd. By rotation and repeated application of $r$, the relation reduces to $X \hookleftarrow S \to \one$ for some subobject $S \hookrightarrow X$. Composing with $X \xleftarrow{\pi_1} X^2 \xrightarrow{\pi_2} X$ (a rotation of $r$) yields $X \xleftarrow{=} X \to \one$. Thus $\cC$ contains $X^n \xleftarrow{=} X^n \to \one$ for each $n$. By coarsening and rotation, we see that $\cC$ contains all relations.
\end{proof}

\begin{proposition} \label{p:r-has-adjoint}
    Consider $\cC = \gtcategory{r}$ and assume that $t \ne 0$. Then $H_\cC = \one \oplus \im(t^{-1}\rho)$, considered as an object of $\cC_\delta$, represents the functor $\cT(-, \one)|_{\cC_\delta}$.
\end{proposition}
\begin{proof}
    It suffices to show that for each $n \in \NN_0$, we have an isomorphism
    \[
        \cC_\delta(X^n, \one) \oplus \cC_\delta(X^n, \im(t^{-1}\rho)) = \cC_\delta(X^n, H_\cC) \longrightarrow \cT(X^n, \one)
    \]
    which on the first summand is given by inclusion and on the second summand by composition with the relation $X \xleftarrow{=} X \to \one$. For $n$ even, we have $\cC_\delta(X^n, H_\cC) = \cC_\delta(X^n, \one) = \cT(X^n, \one)$. For $n = 2m+1$ odd, we need
    \[
        \cC_\delta(X^n, \im(t^{-1}\rho)) \longrightarrow \cT(X^n, \one)
    \]
    to be an isomorphism. For surjectivity, note that every relation $X^n \xhookleftarrow{f} X^b \to \one$ in $\cT(X^n, \one)$ has a preimage $t^{-1} s_f$, where $s_f$ is given by $X^n \xleftarrow{f\pi_1} X^b \times X \xrightarrow{\pi_2} X$.
    For injectivity, we show that $\cC_\delta(X^n, \im(t^{-1}\rho))$ is spanned by relations of the form $s_f$ for subobjects $f$ of $X^n$: Consider a relation $q=(X^n \leftarrow X^c \to X)$ and factor its left leg as $X^c\xrightarrow{\eps}X^b\xhookrightarrow{f}X^n$ with $\eps$ surjective. Composition with $\rho$ gives
    \[\begin{tikzcd}
        && X^c \times_{X} X^2 \arrow{ldd} \arrow{rdd} \arrow[two heads]{d} && \\
        && V \arrow[hook]{d} && \\
        & X^c \arrow{ld} \arrow[pos=.7,swap]{rd} & X^n \times X \arrow[dotted]{lld}\arrow[dotted]{rrd} & X^2 \arrow[pos=.7]{ld} \arrow{rd} & \\
        X^n && X && X \\
    \end{tikzcd} \qquad = \qquad \delta(\eps)s_f.
    \]
    The fiber product at the top is $X^c \times_X X^2 = (X^c \times_X X) \times X$. Its map onto the image $X^b\times X$ is $\eps\times\id_X$, of degree $\delta(\eps)$. Thus $t^{-1}\rho q=t^{-1}\delta(\eps)s_f$, proving the spanning assertion. Injectivity follows as the images $tf$ of the elements $s_f$ are linearly independent.
\end{proof}

\begin{theorem} \label{thm:one-atom} For $\cA$, $\delta$ as above (in particular, $\cA$ has one atom $X$), and $\cC$ a group-theoretical subcategory of $\Rel(\cA)$, assume the following two conditions hold:
\begin{itemize}
    \item[(a)] There is a global splitting object of $\cT_\delta:=\cT(\cA,\delta)$ of the form $\iota(S)$, for $S\in \cC_\delta$, with $\iota\:\cC_\delta\to\cT_\delta$ the inclusion as before.   
    \item[(b)] either \Cref{a:rank-object} holds, or $r\in\cC$.
\end{itemize}
Then $S$ is a global splitting object of $\cC_\delta$, $\cC_\delta$ and $\cT_\delta$ admit abelian envelopes with enough projectives, $\iota$ induces an (exact) tensor functor $\iota\ab\:\cC_\delta\ab\to\cT_\delta\ab$, and $\cC_\delta\ab$ is given by representations of a quantum group in $\cT_\delta\ab$.
\end{theorem}

\begin{proof} If \Cref{a:rank-object} holds, then this follows immediately from \Cref{t:group-theoretical-ab-env} and \Cref{t:group-theoretical-comod}. So let us assume $r\in\cC$ instead. Then by \Cref{l:r-is-special}, we only have to consider two cases: $\cC_1=\Rel(\cA)$ and $\cC_2=\gtcategory{r}$. Let $\cT_1,\cT_2\subset\cT_\delta$ be the corresponding subcategories of the Knop category.
    
Then by \Cref{l:r-is-special}, $\cT_2$ is a subcategory of $\cT_1$ containing all tensor powers of $[X]$ and such that $\Hom_{\cT_2}([X]^{\o a},[X]^{\o b})$ is the full hom-space in $\cT_1$ if $a+b$ is even, or zero if $a+b$ is odd. This implies that any object $A\in\cT_2$ decomposes as a direct sum $A=A_0\oplus A_1$, where $A_i$ consists of those indecomposable summands of $A$ that appear in $[X]^{\o a}$ for $a$ of parity $i$, and that all morphisms between the parity summands of different parity are zero.

We claim that in both cases, $S$ is a global splitting object. This is trivial for $\cT_1$, and follows for $\cT_2$ using the discussed parity decompositions.

We also claim that in both cases, every non-zero morphism $u\:U\to\one$ is a cokernel of $(u\o U-U\o u)$. For $\cT_1$, this follows from \cite{flake2026monoidaladjunctionsabelianenvelopes}*{Proposition~3.10} as $\cT_1=\cT_\delta$ is $\oplus$-pseudo-diagrammatic in the sense of that reference, by \cite{flake2026monoidaladjunctionsabelianenvelopes}*{\S~4.2}. For $\cT_2$, it follows again from the discussed parity decomposition.

With these two claims, the rest of the proofs of \Cref{t:group-theoretical-ab-env} and \Cref{t:group-theoretical-comod} apply. In particular, the central section constructed in the proof of \Cref{t:group-theoretical-comod} works in all cases, and thus \Cref{thm:rec-1} gives the quantum-group assertion, using the right adjoint of $\Ind(\iota\ab)$, which always exists on the level of abelian categories as recalled in \Cref{ss:tannaka-duality}.
\end{proof}

\section{Groups in \texorpdfstring{$\uRep(S_t)$}{Rep(St)}} \label{s:groups-in-rep-st}

One of the first and most well-known examples of Knop's categories $\cT = \cT(\cA, \delta)$ is Deligne's category $\uRep(S_t)$ for $t \in \CC$ a complex parameter (see \cite{deligne2007repst}). It arises by choosing $\cA = \FinSet\op$ to be the opposite category of finite sets, and $\delta(e\op: Y \to X) = t^{\abs{Y \setminus e(X)}}$ for an injective map $e:X\hookrightarrow Y$ of sets. For most of this section, we more concretely take $\cA$ to be the skeleton of $\FinSet\op$, with objects $\underline{n}$ for $n \in \NN_0$. In this case, the calculus of relations $\Rel(\cA)$ is described using partition diagrams. We study the group-theoretical subcategories of $\Rel(\cA)$ and apply the methods from \Cref{s:group-theoretical-subcats} to obtain their abelian envelopes.

However, working with the skeleton of $\FinSet\op$ restricts the class of group-theoretical subcategories. We discuss further choices of the regular category $\cA$ which are equivalent to $\FinSet\op$ but produce new subcategories in \Cref{ss:other-finset-op}.

\subsection{Categories of partitions} \label{ss:categories-of-partitions}
Consider the full skeletal subcategory $\cA$ of $\FinSet\op$ spanned by objects $\underline{n} \coloneqq \{1, \dots, n\}$ for $n \in \NN_0$. Then a relation between $\underline{n} \xleftarrow{f} \underline{b} \xrightarrow{g} \underline{m}$ in $\Rel(\cA)$ is simply a surjection $f \sqcup g: \underline{n+m} \tto \underline{b}$, which we think of as a partition of $n$ upper points and $m$ lower points into $b$ disjoint sets called \emph{blocks}. Vertical and horizontal composition of partition diagrams makes $\Rel(\cA)$ into a rigid strict monoidal category. A more detailed exposition of this diagrammatical calculus is given in \cite{flake2022semisimplicity}. We will use the following notation and terminology:

\begin{itemize}
    \item By a \emph{rotation} of a partition diagram, we mean a repeated application of the isomorphism
    \[
        \Rel(\cA)(\underline{n}, \underline{m}) \cong \Rel(\cA)(\underline{n\pm1}, \underline{m\mp1}).
    \]
    For example,
    \[
        \begin{tikzpicture}[baseline=5pt,scale=0.58]
            \pd{.25,1}; \pd{1,1}; \pd{1.75,1}; \pd{2.5,1};
            \pd{1,0}; \pd{2,0};
            \draw (.25,1) to[out=down,in=down] (1.75,1);
            \draw (1,1) -- (1,0) to[out=up,in=up] (2,0);
            \draw (2.5,.8) -- (2.5,1);
        \end{tikzpicture} \quad\rightsquigarrow\quad
        \begin{tikzpicture}[baseline=5pt,scale=0.58]
            \pd{.25,1}; \pd{1,1}; \pd{1.75,1}; \pd{2.5,1}; \pd{3.25,1};
            \pd{1,0};
            \draw (.25,1) to[out=down,in=down] (1.75,1);
            \draw (1,1) -- (1,0) to[out=up,in=down] (3.25,1);
            \draw (2.5,.8) -- (2.5,1);
        \end{tikzpicture} \quad\rightsquigarrow\quad
        \begin{tikzpicture}[baseline=8pt,scale=0.58]
            \pd{.25,1}; \pd{1,1}; \pd{1.75,1}; \pd{2.5,1}; \pd{3.25,1}; \pd{4,1};
            \draw (.25,1) to[out=down,in=down] (1.75,1);
            \draw (1,1) to[out=down,in=down] (3.25,1) to[out=down,in=down] (4,1);
            \draw (2.5,.8) -- (2.5,1);
        \end{tikzpicture}
    \]
    \item By \emph{applying a cap} to a partition diagram, we mean composition with the diagram \tikz[baseline=(origin.base),scale=0.5]{
        \coordinate (origin) at (0,0); 
        \pd{0,0}; \pd{1,0}; \pd{1.5,0}; \pd{2,0}; \pd{2.5,0}; \pd{3.5,0};
        \pd{0,.5}; \pd{1,.5}; \pd{2.5,.5}; \pd{3.5,.5};
        \draw (0,0) -- (0,.5); \draw (1,0) -- (1,.5); \draw (2.5,0) -- (2.5,.5); \draw (3.5,0) -- (3.5,.5);
        \draw (1.5,0) to[out=up,in=up] (2,0);
        \node (d1) at (.5,.25) {\tiny $\dots$};
        \node (d2) at (3,.25) {\tiny $\dots$};
    }, where the cap is placed in an appropriate position.
    \item By a \emph{category of partitions}, we mean a rigid monoidal subcategory $\cC$ of $\Rel(\cA)$ containing all objects $\underline{n}$, $n \in \NN_0$ which is closed under horizontally flipping partition diagrams.
    \item Let $P$ be a set of partition diagrams. We denote by $\langle{}P\rangle{}$ the smallest category of partitions containing $P$.
    \item We call a category of partitions \emph{group-theoretical} if it is closed under coarsening (i.e.~joining blocks in partitions). 
\end{itemize}

A subcategory of $\Rel(\cA)$ is a group-theoretical category of partitions if and only if it is a group-theoretical category of relations in the sense of \Cref{d:group-theoretical}. Any category of partitions gives rise to a (non-full) subcategory $\uRep(\cC, t)$ of $\uRep(S_t)$. For $\cC$ group-theoretical, this is precisely the category $\cC_{\delta_t}$ from \Cref{d:c-delta}. We denote the image of $\underline{n}$ in $\uRep(\cC, t)$ by $[n]$.

The construction of $\uRep(\cC, t)$ comes from the construction and classification of easy quantum groups initiated by Banica and Speicher \cite{banica2009liberation} and completed by Raum and Weber (see \cite{raum2016full} and references therein). The term \emph{group-theoretical} reflects the fact that such categories of partitions are determined by certain subgroups of the infinite free product $\ZZ_2^{* \infty}$. The classification implies that there are uncountably many categories of partitions, and all but countably many are group-theoretical.

\begin{lemma}[\cite{raum2014combinatorics}*{Lemma~2.3}]
    A category of partitions is group-theoretical if and only if it contains the partition \coarsener. In other words, $\gtcategory{P} = \langle{}P, \coarsener\rangle{}$ for any set of partitions $P$.
\end{lemma}

Recall from \Cref{ss:one-atom} that the partition $r = \upperdots{2}$ plays a special role: If a group-theoretical category $\cC$ of partitions contains $\upperdots{2}$, then $\cC = \Rel(\cA)$ or $\cC = \langle{}\upperdots{2}, \coarsener\rangle{}$, and for $t\ne0$ in these cases $\iota: \uRep(\cC, t) \to \uRep(S_t)$ has a right adjoint. We can classify categories of partitions containing $\crossing$:

\begin{lemma} \label{l:crossing-partitions}
    If a group-theoretical category $\cC$ of partitions contains the crossing $\crossing$, then $\cC = \Rel(\cA)$ or $\cC = \langle{}\upperdots{2}, \coarsener\rangle{}$ or $\cC = \langle{}\crossing, \coarsener\rangle{}$.
\end{lemma}
\begin{proof}
    We have a chain of (non-full) subcategories
    \[\begin{tikzcd}
        \cC_0 = \langle{}\crossing, \coarsener\rangle{} \arrow[hook]{r} &
        \cC_1 = \langle{}\upperdots{2}, \coarsener\rangle{} \arrow[hook]{r} &
        \Rel(\cA).
    \end{tikzcd}\]
    Indeed, $\crossing$ arises from $\upperdots{4}$ by coarsening and rotating down two of the points, so $\crossing \in \cC_1$. It suffices to show that $\cC_1$ is the only category of partitions sitting between $\cC_0$ and $\Rel(\cA)$.

    Consider a category of partitions $\cC$ containing $\cC_0$, let $p$ be a partition in $\cC$ and let $n$ be the number of blocks in $p$ of odd size. By rotation, we may assume that it has only upper points. By permuting the points, we may assume that points can only be connected if they are adjacent. By applying caps to blocks of size $\geq 2$, we may assume that $p$ is the partition consisting of $n$ discrete dots. For example,
    \[
        \begin{tikzpicture}[baseline=5pt,scale=0.58]
            \pd{.25,1}; \pd{1,1}; \pd{1.75,1}; \pd{2.5,1};
            \pd{1,0}; \pd{2,0};
            \draw (.25,1) to[out=down,in=down] (1.75,1);
            \draw (1,1) -- (1,0) to[out=up,in=up] (2,0);
            \draw (2.5,.8) -- (2.5,1);
        \end{tikzpicture} \in \cC \xRightarrow{\text{rotate}}
        \begin{tikzpicture}[baseline=8pt,scale=0.58]
            \pd{.25,1}; \pd{1,1}; \pd{1.75,1}; \pd{2.5,1}; \pd{3.25,1}; \pd{4,1};
            \draw (.25,1) to[out=down,in=down] (1.75,1);
            \draw (1,1) to[out=down,in=down] (3.25,1) to[out=down,in=down] (4,1);
            \draw (2.5,.8) -- (2.5,1);
        \end{tikzpicture} \in \cC \xRightarrow{\text{permute}}
        \begin{tikzpicture}[baseline=8pt,scale=0.58]
            \pd{.25,1}; \pd{1,1}; \pd{1.75,1}; \pd{2.5,1}; \pd{3.25,1}; \pd{4,1};
            \draw (.25,1) to[out=down,in=down] (1,1);
            \draw (2.5,1) to[out=down,in=down] (3.25,1) to[out=down,in=down] (4,1);
            \draw (1.75,.8) -- (1.75,1);
        \end{tikzpicture} \in \cC \xRightarrow{\text{apply caps}}
        \begin{tikzpicture}[baseline=8pt,scale=0.58]
            \pd{.25,1}; \pd{.75,1};
            \draw (.25,.8) -- (.25,1);
            \draw (.75,.8) -- (.75,1);
        \end{tikzpicture} \in \cC.
    \]
    Again applying caps, we see that $\upperdots{2} \in \cC$ if $n \geq 2$ is even and $\upperdots{1} \in \cC$ if $n$ is odd. In total, we can distinguish three cases:

    \begin{itemize}
        \item The category $\cC_0$ contains all Brauer diagrams and their coarsenings. Thus it consists precisely of the diagrams with even block size. If all diagrams in $\cC$ have $n = 0$ odd blocks, then $\cC = \cC_0$.
        \item The category $\cC_1$ contains all diagrams with an even total of (upper and lower) points. Thus it consists precisely of the diagrams with an even number $n$ of odd blocks. If all diagrams in $\cC$ have $n$ even but some have $n > 0$, then $\cC = \cC_1$.
        \item Finally, if some diagram in $\cC$ has $n$ odd, then $\cC = \Rel(\cA)$. \qedhere
    \end{itemize}
\end{proof}

\begin{lemma} \label{l:no-single-points}
    Let $\cC$ be a group-theoretical category of partitions. If $\cC$ contains a partition $p$ with a block of size 1, then $\upperdots{2} \in \cC$.
\end{lemma}
\begin{proof}
    By rotation, we may assume that $p \in \cC(\underline{n}, \underline{1})$ for some $n$ and the single bottom point is not connected to any other point. Applying caps to the top points, we see that $\upperdots{2} \in \cC$ if $n$ is odd and $\upperdots{1} \in \cC$ if $n$ is even.
\end{proof}

\subsection{Existence of abelian envelopes} \label{ss:existence-of-abelian-envelopes}

\begin{proposition}[{see \cite{flake2022semisimplicity}}]
    Let $\cC$ be a group-theoretical category of partitions. Then $\uRep(\cC, t)$ is semisimple if and only if $t \notin \NN_0$.
\end{proposition}

For $\cC = \Rel(\cA)$, this generalizes a fundamental property of $\uRep(S_t)$ (see~\cite{deligne2007repst}*{Théorème~2.18}). For $d \in \NN_0$, $\uRep(S_d)$ is not abelian but admits an abelian envelope $\uRep\ab(S_d)$ (see \cite{comes2014repabst}). To see the existence of an abelian envelope for $\uRep(\cC, d)$ for any group-theoretical category $\cC$, we check \Cref{a:rank-object}.

\begin{remark} \label{rem:upperdots}
 Recall that by \Cref{l:r-is-special}~(b), there are only two group-theoretical categories of partitions containing $\upperdots{2}$: $\cC = \Rel(\cA)$ and $\cC = \gtcategory{\upperdots{2}}$.
\end{remark}

\begin{lemma} \label{l:no-two-upperdots}
    Let $\cC$ be a group-theoretical category of partitions other than $\gtcategory{\upperdots{2}}$. Then $\cC$ satisfies \Cref{a:rank-object}.
\end{lemma}
\begin{proof}
    We assume first $\upperdots{2}\not\in\cC$. Consider a relation $\underline{n} \xleftarrow{f} \underline{b} \xrightarrow{g} \underline{m}$ in $\Rel(\cA)$ and its corresponding partition diagram. Then $f$ is injective in $\cA$ if and only if there are no bottom blocks (i.e.~blocks without upper points) in the partition diagram.

    Fix an injection $f: \underline{n} \hookrightarrow \underline{m}$ (i.e. a surjective map of sets) and define $R_f \coloneqq \underline{r_f}$, where
    \[
        r_f \coloneqq \min\{i \in \NN_0 \mid \exists \text{relation in $\cC$ of the form } \underline{m} \xhookleftarrow{f} \underline{n} \twoheadleftarrow \underline{b} \to \underline{i}\}.
    \]
    This minimum is well-defined: The set is non-empty because the relation $\beta_f$ given by $\underline{m} \xhookleftarrow{f} \underline{n} \xhookrightarrow{f} \underline{m}$ is a coarsening of $\id_{\underline{m}}$ and therefore contained in $\cC$.

    We need to show that $R_f$ satisfies the properties (A1) and (A2) in \Cref{a:rank-object}. Consider a relation in $\cC$ of the form $\underline{m} \xhookleftarrow{f} \underline{n} \stackrel{\eps}{\twoheadleftarrow} \underline{b} \to \underline{i}$. We show that if $\eps$ is not a bijection, then $i > r_f$. If $\eps$ is not a bijection, then the corresponding diagram has a bottom block. By \Cref{l:no-single-points}, we have $i \geq 2$. Therefore, we can choose a point in a bottom block and one of its neighbors and apply a cup to them. In pictures:
    \[
        \begin{tikzpicture}[baseline=5pt,scale=0.58]
            \draw[gray] (.25,0) to[out=down,in=down] (1,0);
            \pd{.25,1}; \pd{1,1}; \pd{1.75,1};
            \pd{.25,0}; \pd{1,0}; \pd{1.75,0};
            \draw (1,1) to[out=down,in=down] (1.75,1);
            \draw (.25,1) to[out=down,in=up] (1,0);
            \draw (.25,0) to[out=up,in=up] (1.75,0);
        \end{tikzpicture} \in \cC \xRightarrow{\text{apply cup to bottom block}}
        \begin{tikzpicture}[baseline=5pt,scale=0.58]
            \pd{.25,1}; \pd{1,1}; \pd{1.75,1};
            \pd{1,0};
            \draw (1,1) to[out=down,in=down] (1.75,1);
            \draw (.25,1) to[out=down,in=up] (1,0);
        \end{tikzpicture} \in \cC
    \]
    This produces a partition in $\cC(\underline{m}, \underline{i-2})$ which has the same top connectivity as the original partition but fewer bottom points.

    This shows that $R_f$ satisfies (A2). Together with the definition of $R_f$, this implies (A1), under the assumption $\upperdots{2}\not\in\cC$.

    For $\cC=\Rel(\cA)$, defining $R_f \coloneqq \underline{0}$ works for all $f$.

    \Cref{rem:upperdots} implies that the special cases discussed actually imply the assertion.
\end{proof}

In~\cite{comes2014repabst}*{Lemma~3.11} it was shown that the object $\Delta_{d+1} \coloneqq \im(x(\id_{[d+1]}))$ splits all morphisms in $\uRep(S_d)$. In terms of the (oligomorphic) infinite symmetric group $G$ acting on $\Omega = \{1, 2, \dots\}$, the object $\Delta_{d+1}$ corresponds to the $G$-set $G/U$ where $U$ is the subgroup fixing $1, \dots, d+1$ pointwise. It is contained in all $\cC_\delta$ for $\cC \subseteq \Rel(\cA)$ group-theoretical, as its defining idempotent is a linear combination of coarsenings of an identity diagram.

\begin{corollary} \label{c:partitions-ab-env}
    Let $\cC$ be a group-theoretical category of partitions. For $d \in \NN_0$, $\Delta_{d+1}$ is a global splitting object already in $\uRep(\cC, d)$, and $\uRep(\cC, d)$ admits an abelian envelope with enough projectives. This abelian envelope is given by representations of a quantum group in $\uRep(S_d)\ab$.
\end{corollary}
\begin{proof} This is \Cref{thm:one-atom} using \Cref{l:no-two-upperdots}.
\end{proof}

\subsection{Interpretation as quantum groups} \label{ss:quantum-groups}

Let $\cC$ be a group-theoretical category of partitions. \Cref{c:partitions-ab-env} yields an abelian envelope $\uRep(\cC, t)\ab$. The construction of the central section in the proof of \Cref{t:group-theoretical-comod} does not require \Cref{a:rank-object}. Thus \Cref{thm:rec-1} realizes $\uRep(\cC, t)\ab$ as $\Comod_{\uRep(S_t)\ab}(H)$ for some Hopf algebra Ind-object in $\uRep(S_t)\ab$, or alternatively as a ``quantum group'' in $\uRep(S_t)\ab$.

If $\uRep(\cC, t)$ and $\iota: \uRep(\cC, t) \to \uRep(S_t)$ are symmetric, $H$ becomes a commutative Hopf algebra and can therefore be seen as the coordinate ring of an affine group over $\uRep(S_t)\ab$. For symmetry, we need $\crossing \in \cC$, which by \Cref{l:crossing-partitions} means that $\cC$ is one of the following categories of partitions:

\begin{itemize}
    \item For $\cC = \Rel(\cA)$, we have $H_\cC = \one$, which corresponds to the trivial group.
    \item For $\cC = \gtcategory{\upperdots{2}}$, and $t \ne 0$, we have $H = \iota(H_\cC) = \one \oplus \im(t^{-1} \verticaldots) \cong \one \oplus \one \cong \cO(S_2)$ in $\uRep(S_t)$. Thus $\uRep(\cC, t)\ab$ corresponds to the discrete group $S_2$ defined over $\uRep(S_t)\ab$. Note that the categories $\uRep(\cC, t)$ for general $t$ interpolate the representation categories $\Rep(S_2 \times S_n) = \Rep_{\Rep(S_n)}(S_2)$ for $n \geq 1$.
    \item For $\cC = \gtcategory{\crossing}$, we have $H = \bigoplus_{n \in \NN_0} \im(e_n)$, where $e_n$ is the idempotent endomorphism of $[n]$ given by
    \begin{equation} \label{e:idempotent-hyperoctahedral}
        e_n \coloneqq \left(\frac{1}{n!} \sum_{\sigma \in S_n} \sigma\right) x(\id_{[n]}).
    \end{equation}
    In the following, we show that the Hopf algebra $H$ is isomorphic to $\kk[S_2]^{\o t}$, an object which after semisimplification at $t=d \in \NN_0$ becomes the group algebra $\kk[S_2]^{\o d} = \kk[S_2^d]$, or equivalently the algebra of functions on $S_2^d$. Thus $\uRep(\cC, t)\ab$ can be seen as ``representations of the group $S_2^t$ over $\uRep(S_t)\ab$''.
\end{itemize}

In \cite{reptheocomplexrank}, a symmetric monoidal functor $(-)^{\otimes t}: \mathrm{UVect} \to \Ind(\uRep(S_t))$ was constructed for $t\notin\NN_0$, where $\mathrm{UVect}$ is the category of $\CC$-vector spaces with a distinguished non-zero vector and linear maps preserving this vector. It interpolates the functors $(-)^{\otimes d}: \mathrm{UVect} \to \Rep(S_d)$ for $d \in \NN_0$. Being symmetric monoidal, it sends algebras to algebras and Hopf algebras to Hopf algebras. Note that this functor is not linear.

\begin{proposition}[{see \cite{reptheocomplexrank}*{Proposition~4.3}}] \label{p:gens-rels-tth-power}
    For a unital $\kk$-algebra $A$, consider the tensor algebra $T_A \coloneqq T(A \o [1])$ with graded pieces $T_A^i \coloneqq (A \o [1])^{\o i}$. Then
    the algebra $A^{\o t}$ is the quotient of $T_A$ by the ideal generated by the images of the morphisms
    \[\begin{aligned}
        &r_1 = \eta_A\o\id_{[1]}-\upperdots{1}: [1]\to T_A^1\oplus T_A^0, \\
        &r_2 = \mu_A \o \id_{[1]} - \sigma_{23}(\id_{A \o A} \o \;\comult) : A \o A \o [1] \to T_A^1 \oplus T_A^2, \\
        &r_3 = \id - \sigma_{13}\sigma_{24} - ([-,-]_A \o \mult)\sigma_{23} : A \o [1] \o A \o [1] \to T_A^1 \oplus T_A^2.
    \end{aligned}\]
    Here $\eta_A$, $\mu_A$ and $[,]_A$ denote the unit, multiplication and commutator of $A$, respectively, and $\sigma_{ij}$ permutes the $i$-th and $j$-th tensor factor. The term $\upperdots{1}$ in $r_1$ takes values in $T_A^0=\one$; it identifies the unit of $A$ with the unit of the tensor algebra.
\end{proposition}

For $t\in\NN_0$, we define $A^{\o t}$ in $\Ind(\uRep(S_t)\ab)$ by this presentation.

\begin{lemma} \label{l:gens-rels-ks2-otimes-t}
    Consider $A = \kk[S_2]$. The Hopf algebra $A^{\o t}$ is isomorphic to the quotient of $T([1])$ by the ideal generated by the images of
    \[
        s_1 = \id_{[2]} - \;\crossing: [2] \to [2], \qquad s_2 = \comult - \upperdots{1}: [1] \to [2] \oplus [0].
    \]
    The comultiplication and counit are induced by $\comult$ and $\upperdots{1}$, respectively.
\end{lemma}
\begin{proof}
    Write $A \coloneqq \kk[S_2] = \kk[z]/(z^2-1)$ and $A \o [1] = 1 \o [1] \oplus z \o [1]$. We use the presentation of $A^{\o t}$ from \Cref{p:gens-rels-tth-power}. We set $Q=T([1])/(s_1,s_2)$ and let $j:[1]\to Q$ be the canonical morphism and $\eta_Q\:\one\to Q$ the unit. The morphism $A\o[1]\to Q$ given by $\eta_Q\upperdots{1}$ on $1\o[1]$ and by $j$ on $z\o[1]$ induces an algebra morphism $\sigma:A^{\o t}\to Q$: relation $r_1$ holds by definition, $r_2$ follows from $s_2$ and the counit identities, and $r_3$ follows from $s_1$.

    Conversely, the inclusion $z\o[1]\subseteq A\o[1]$ induces an algebra morphism $\tau:Q\to A^{\o t}$. Indeed, $r_3$ gives $s_1$, while $r_2$ on $z\o z\o[1]$, together with $r_1$, gives $s_2$. The relations imply that the compositions $\sigma\tau$ and $\tau\sigma$ are identities on the generators.    
    One checks on the generator $z\o[1]$ that $\tau$ preserves comultiplication and counit and is therefore an isomorphism of Hopf algebras.
\end{proof}

\begin{proposition} \label{p:hopf-alg-iso}
    Consider the category $\cC = \gtcategory{\crossing}$ of partitions with blocks of even size. Then $\uRep(\cC, t)\ab \cong \Comod_{\uRep(S_t)\ab}(H)$, where $H$ is the Hopf algebra $\kk[S_2]^{\otimes t} \in \Ind(\uRep(S_t)\ab)$.
\end{proposition}
\begin{proof}
    Abbreviate $\cC_t \coloneqq \uRep(\cC, t)$ and $\cS_t \coloneqq \uRep(S_t)$. We show that the Hopf algebra $T([1])/(s_1, s_2)$ from \Cref{l:gens-rels-ks2-otimes-t} is isomorphic to the Hopf algebra $T(W)/\ker(\psi)$ from \Cref{thm:rec-2}, where we choose $W'=[1]\in\cC_t$ and $W = [1]$ with the coalgebra structure $\Delta = \comult$ and $\eps = \upperdots{1}$. The morphism $\Delta':W'\to W'\o J(W)$ is given by the single block on one upper and three lower points. Under the identification $\iota J(W)\cong W$, its image is $\Delta$. To describe $\psi$ explicitly, recall that the right adjoint of $\Ind(\iota: \cC_t \to \cS_t)$ is built from the isomorphism
    \begin{equation} \label{e:adjunction-iso}
        \Ind(\cC_t)\left(-, \bigoplus_{n \geq 0} \im(e_n)\right) \cong \Ind(\cS_t)\left(\Ind(\iota)(-), [0]\right)
    \end{equation}
    from \Cref{p:representing-object}. Here $e_n$ is the idempotent given in \Cref{e:idempotent-hyperoctahedral}. For fixed $k \geq 0$, consider a partition of the form $f: [k] \to [0]$ with $o(f)$ odd blocks. The partition $s_f: [k] \to [o(f)]$ connects every odd block of $f$ with a new bottom point. Under \Cref{e:adjunction-iso}, the Möbius basis element $x(f): [k] \to [0]$ has preimage $e_{o(f)}x(s_f)$. Möbius inversion states that $\eps^{\o k} = \sum_{f: [k] \to [0]} x(f)$, so the adjunction mate of $\eps^{\o k}$ under \Cref{e:adjunction-iso} is
    \[
        \eps^{\o k,\bot} = \sum_{f: [k] \to [0]} \left(e_{o(f)}x(s_f): [k] \to \im(e_{o(f)})\right).
    \]
    On the graded piece $[k]$, the morphism $\psi: T([1]) \to \bigoplus_{n \geqslant 0} \im(e_n)$ is given by $\eps^{\o k,\bot}$. In particular, we have $\eps^{\o 0,\bot} = \id_{\im(e_0)} = \id_{[0]}$ and $\eps^{\o 2,\bot} = \cuppartition + e_2: [2] \to [0] \oplus \im(e_2)$. Thus
    \[
        \psi \circ s_1 = \eps^{\o 2,\bot} - \eps^{\o 2,\bot} \circ \crossing = \cuppartition - \tikz[baseline=(origin.base),scale=0.5]{
            \coordinate (origin) at (0,0);
            \pd{0,.5}; \pd{0,0}; \pd{.5,0}; \pd{.5,.5};
            \draw (0,.5) -- (.5, 0) to[out=down,in=down] (0,0) -- (.5, .5);
        } + e_2 - e_2 \circ \crossing = 0
    \]
    and
    \[
        \psi \circ s_2 = \eps^{\o 2,\bot} \circ \comult - \eps^{\o 0,\bot} \circ \upperdots{1} = \tikz[baseline=(origin.base),scale=0.5]{
            \coordinate (origin) at (0,0);
            \pd{.25,.5}; \pd{0,0}; \pd{.5,0};
            \draw (.25,.5) -- (.25, .25) -- (0,0) to[out=down,in=down] (.5,0) -- (.25, .25);
        } + e_2 \circ \comult - \upperdots{1} = 0.
    \]
    In total, $\psi$ factors as
    \[\begin{tikzcd}
        T([1]) \arrow{rr}{\psi} \arrow[two heads]{rd}{q} && \bigoplus_{n \geq 0} \im(e_n). \\
        & T([1])/(s_1, s_2) \arrow{ru}{\overline{\psi}} &
    \end{tikzcd}\]
    A section of $\psi$ is given by the morphism $\vfi: \bigoplus_n \im(e_n) \to T([1])$ including each $\im(e_n) \subseteq [n]$ as a direct summand. Indeed, in the composition $\eps^{\o n,\bot} e_n$, every summand is killed except the finest partition belonging to $\upperdots{1}^{\o n}$, so $\eps^{\o n,\bot} e_n = \id_{\im(e_n)}$. It remains to show that $q\vfi\psi=q$.
    
    For fixed $k$ and a partition $f: [k] \to [0]$, consider the morphism $q x(\beta_f)$. By applying the symmetry relation $s_1$, we can reorder the bottom blocks of $\beta_f$ in such a way that connected bottom points are adjacent. Next, repeatedly apply relation $s_2$ to each block, deleting connected bottom points in pairs. This leaves one bottom point for each odd block of $f$ and no bottom points for each even block of $f$. We see that $q x(\beta_f) = q \vfi e_{o(f)} x(s_f)$. Thus
    \[
        q \id_{[k]} = q \sum_{f: [k] \to [0]} x(\beta_f) = q \vfi \sum_{f: [k] \to [0]} e_{o(f)} x(s_f) = q \vfi \eps^{\o k,\bot},
    \]
    so $q = q \vfi \psi=q\vfi\overline{\psi}q$. Because $q$ is an epimorphism, $q\vfi\overline{\psi}=\id$. Together with $\overline{\psi}q\vfi=\psi\vfi=\id$, this shows that $\overline{\psi}$ is an isomorphism with inverse $q\vfi$.
\end{proof}

\subsection{Other realizations of \texorpdfstring{$\FinSet\op$}{(FinSet)op}} \label{ss:other-finset-op} We can replace $\cA$ by a different category that is still equivalent to $\FinSet\op$. Concretely, consider a word $\mathbf{v} = v_1 \dots v_n$ over a non-empty alphabet $A$ and the set $[\mathbf{v}] \coloneqq \{(1, v_1), \dots, (n, v_n)\}$. Let $\cB$ be the full subcategory of $\FinSet\op$ spanned by the sets $[\mathbf{v}]$ where $\mathbf{v}$ runs over all words in $A^*$, with chosen direct products $[\mathbf{v}] \times [\mathbf{w}] = [\mathbf{v}\mathbf{w}]$ (concatenation of words).

A relation $[\mathbf{v}] \xleftarrow{f} [\mathbf{b}] \xrightarrow{g} [\mathbf{w}]$ in $\Rel(\cB)$ is a quotient of sets $f \sqcup g: [\mathbf{vw}] \tto [\mathbf{b}]$. It does not depend on the word $\mathbf{b}$ but only on its length. In other words, a relation can be depicted as a partition diagram whose upper and lower points are labelled by elements of the alphabet $A$. For $A$ a singleton set, we recover the unlabelled partitions $\Rel(\cA)$. As a category, $\cB$ is equivalent to $\cA$ and $\FinSet\op$, but we find new group-theoretical subcategories in $\Rel(\cB)$ compared to $\Rel(\cA)$.

\begin{example} \label{ex:klein-four-group}
    Assume that the alphabet $A$ consists of two elements $\{a, b\}$ and consider the subcategory $\cC = \gtcategory{\tikz[scale=0.5]{
        \pdg{0,.5}{south}{a};
        \draw (0,.3) -- (0,.5);
        \pdg{0.5,.5}{south}{a};
        \draw (0.5,.3) -- (0.5,.5);
    },  \tikz[scale=0.5]{
        \pdg{0,.5}{south}{b};
        \draw (0,.3) -- (0,.5);
        \pdg{0.5,.5}{south}{b};
        \draw (0.5,.3) -- (0.5,.5);
    }, \tikz[scale=0.5,baseline=0pt]{
        \pdg{0,.5}{south}{a}; \pdg{.5,.5}{south}{b};
        \pdg{0,0}{north}{b}; \pdg{.5,0}{north}{a};
        \draw (0,.5) -- (.5,0); \draw (.5,.5) -- (0,0);
    }}$ of $\Rel(\cB)$. The same-colour crossings follow from the two-singleton generators by coarsening and rotation, so $\cC$ is symmetric. Give each word its colour-parity degree in $\Gamma=\ZZ_2\times\ZZ_2$. Then $\cC$ consists precisely of all partitions whose upper and lower words have the same degree: all generators and closure operations preserve this degree, and every such partition is obtained by reordering and rotating products of the two-singleton generators and then coarsening.

    For $t \ne 0$, the functor $\iota: \cC_t \to \uRep(S_t)$ has a right adjoint sending $\one$ to the object
    \[
        H_\cC = \one \oplus \im\left(t^{-1} \tikz[baseline=(origin.base),scale=0.5]{
            \coordinate (origin) at (0,0);
            \pdg{0,0}{north}{a}; \pdg{0,.5}{south}{a};
            \draw (0,0) -- (0,.18); \draw (0,.32) -- (0,.5);
        }\right) \oplus \im\left(t^{-1} \tikz[baseline=(origin.base),scale=0.5]{
            \coordinate (origin) at (0,0);
            \pdg{0,0}{north}{b}; \pdg{0,.5}{south}{b};
            \draw (0,0) -- (0,.18); \draw (0,.32) -- (0,.5);
        }\right) \oplus \im\left(t^{-1} \tikz[baseline=(origin.base),scale=0.5]{
            \coordinate (origin) at (0,0);
            \pdg{0,0}{north}{a}; \pdg{0,.5}{south}{a};
            \pdg{.5,0}{north}{b}; \pdg{.5,.5}{south}{b};
        \draw (0,.5) to[out=down,in=down] (.5,.5);
        \draw (0,0) to[out=up,in=up] (.5,0);
        }\right)
    \]
    Indeed, the four summands $E_g$, indexed by $g\in\Gamma$, have distinct degrees and satisfy $\iota(E_g)\cong\one$. Morphism spaces within each degree agree with those in $\uRep(S_t)$, whereas morphisms between different degrees vanish. Hence
    \[
        \cC_t(-,H_\cC)\cong\uRep(S_t)(\iota(-),\one).
    \]
    Rigidity and dominance then give the asserted right adjoint.

    The summands satisfy $E_g\o E_h\cong E_{g+h}$, and $\cC_t$ identifies with the category of $\Gamma$-graded objects of $\uRep(S_t)$, with $\iota$ forgetting the grading. Thus $\iota(H_\cC)\cong\kk[\Gamma]\cong\cO(\Gamma)$ is the coordinate algebra of the Klein four-group over $\uRep(S_t)\ab$. Such a group does not arise from group-theoretical subcategories of $\Rel(\cA)$.
\end{example}

\section{Groups in \texorpdfstring{$\uRep(G \wr S_t)$}{Rep(G wr St)}} \label{s:groups-in-rep-g-st}

The statements in \Cref{s:groups-in-rep-st} admit a slight generalization: Instead of $\FinSet\op$, consider the opposite of the category of free finite-rank $G$-sets for a fixed finite group $G$. Its skeleton $\cA$ consists of finite coproducts $G^{\sqcup n}$ for $n \in \NN_0$. Writing $\ell_0(G^{\sqcup n}) \coloneqq n$, all degree functions on $\cA$ are of the form $\delta(e\op: Y \to X) = t^{\ell_0(Y \setminus e(X))}$ for injective maps $e:X\hookrightarrow Y$ of free finite $G$-sets and some fixed $t \in \kk$. The categories $\cT(\cA, \delta)$ interpolate the categories of representations of wreath products $\Rep(G \wr S_n)$ \cite{knop}*{Example~2 in \S~9} and are therefore also denoted by $\uRep(G \wr S_t)$. They are semisimple if and only if $t \notin \NN_0 \abs{G}$ \cite{knop}*{Example~2 in \S~8}.

The calculus of relations $\Rel(\cA)$ has been described in terms of \emph{$G$-partitions} in~\cite{savage2021group}*{\S~3}. A $G$-partition diagram is a partition diagram (its \emph{shape}) with top and bottom points
labeled by elements of $G$. Two $G$-partition diagrams are equivalent if
for each block, the labels of one diagram can be obtained from the labels
of the other by multiplying with a fixed group element. Relations in
$\Rel(\cA)$ and their composition correspond precisely to equivalence
classes of $G$-partition diagrams and their composition.

A group-theoretical subcategory $\cC \subseteq \Rel(\cA)$ contains all
$G$-partition diagrams of the form
\[\begin{tikzpicture}[baseline=2pt]
    \pdg{0,.3}{south}{a};
    \pdg{0,0}{north}{b};
    \draw (0,0) -- (0,.3);
\end{tikzpicture}
\qquad \text{for all }a, b \in G.
\]
Consequently, $\cC$ is entirely determined by the shapes of the
$G$-partition diagrams it contains, i.e.~if a $G$-partition diagram is in
$\cC$, then all $G$-partition diagrams with the same shape are also in
$\cC$. Therefore, the same arguments as in \Cref{ss:categories-of-partitions,ss:existence-of-abelian-envelopes} show that $\cC$ satisfies \Cref{a:rank-object} if it does not contain the two-singleton relation. Using the following input from \cite{flake2026abelianenvelopesinterpolationcategories}*{\S~4.2}, we can deduce that $\cC_\delta$ has an abelian envelope.

For $t \in \kk$ and $n \in \NN_0$, consider the elements
\[
    x_n \coloneqq x(\id_{G^{\sqcup n}}), \qquad e_n \coloneqq \frac{1}{n!} \sum_{\sigma \in S_n} \sigma, \qquad \omega_n \coloneqq \frac{1}{\abs{G}^n} \sum_{g_1, \dots, g_n \in G} \begin{tikzpicture}[baseline=2pt]
        \pdg{0,.3}{south}{1};
        \pdg{0,0}{north}{g_1};
        \draw (0,0) -- (0,.3);
        \node at (.4, .15) {\dots};
        \pdg{.8,.3}{south}{1};
        \pdg{.8,0}{north}{g_n};
        \draw (.8,0) -- (.8,.3);
    \end{tikzpicture}.
\]
Then $e \coloneqq x_ne_n\omega_n$ is an idempotent in $\uRep(G \wr S_t)$ and we denote $\Delta_n \coloneqq \im(e)$. There exists a monoidal ``restriction'' functor
\[
    \uRep(G \wr S_t) \longrightarrow \uRep(G \wr S_{t-n \abs{G}}) \boxtimes \Rep(G \wr S_n),
\]
which has a two-sided adjoint sending $\one \mapsto \Delta_n$.

Now consider $d \in \NN_0 \abs{G}$ and $n = d+1$. Then $\uRep(G \wr S_d)$ is not abelian but $\uRep(G \wr S_{d-n\abs{G}}) \boxtimes \Rep(G \wr S_n)$ is even semisimple. By \Cref{l:adjoint-preserves-splitting-objects}, this implies that $\Delta_n = \Delta_{d+1}$ is a global splitting object in $\uRep(G \wr S_d)$, which therefore admits an abelian envelope with enough projectives. For any group-theoretical subcategory $\cC \subseteq \Rel(\cA)$, the functor $J:\cT\to\cC_\delta$ constructed in the proof of \Cref{t:group-theoretical-comod} gives a lift $J(\Delta_n)$ with $\iota J(\Delta_n)\cong\Delta_n$. Its construction does not require abelian envelopes, and the symmetrizer $e_n$ itself need not belong to $\cC_\delta$.

\begin{corollary} \label{c:wreath-products-ab-env}
    Let $d \in \NN_0 \abs{G}$ with associated degree function $\delta$. For any group-theoretical subcategory $\cC \subseteq \Rel(\cA)$, the pseudo-tensor category $\cC_\delta$ has a global splitting object whose image under $\iota$ is isomorphic to $\Delta_{d+1}$ in $\cT$. It admits an abelian envelope with enough projectives which is given by representations of a quantum group in $\uRep(G\wr S_d)\ab$.
\end{corollary}
\begin{proof}
    This follows from \Cref{thm:one-atom} with $S=J(\Delta_{d+1})$.
\end{proof}

\begin{remark}
    Instead of free $G$-sets, we can consider all finite $G$-sets. The resulting Knop category and its interpolation property were examined in \cite{entovaaizenbud2024deligneknoptensorcategoriesfunctoriality}*{\S~3.5.1}. In contrast to $\FinSet\op$ and $G\mathrm{-FreeFinSet}\op$ where all objects are direct products of a single atomic object (\Cref{ss:one-atom}), the category $G\mathrm{-FinSet}\op$ has an atom for each isomorphism class of transitive $G$-sets. We consider the following realization of $G\mathrm{FinSet}\op$ which allows for a diagrammatic calculus.
    
    Let $A = \{T_1, \dots, T_r\}$ be a full set of representatives of isomorphism classes of transitive $G$-sets. For each word $\mathbf{v} = T_{i_1} \cdots T_{i_n}$ over the alphabet $A$, define the $G$-set 
    \[
        [\mathbf{v}] \coloneqq \{(t, j) \mid j \in \{1, \dots, n\}, t \in T_{i_j}\} \cong \bigsqcup_{j=1}^n T_{i_j}
    \]
    and let $\cA$ be the full subcategory of $G\text{-}\mathrm{Set}\op$ spanned by such $[\mathbf{v}]$. This is a regular category with products given by $[\mathbf{v}] \times [\mathbf{w}] \coloneqq [\mathbf{v}\mathbf{w}]$ (concatenation of words). For words $\mathbf{v} = T_{i_1} \cdots T_{i_n}, \mathbf{w} = T_{\ell_1} \cdots T_{\ell_m}$, a morphism $[\mathbf{v}] \to [\mathbf{w}]$ in $\cA$ is uniquely determined by a map of sets $f: \{1, \dots, m\} \to \{1, \dots, n\}$ and morphisms of $G$-sets $\vfi_j: T_{\ell_j} \to T_{i_{f(j)}}$ which are necessarily surjective. A relation in $\Rel(\cA)$ may be represented by a partition diagram with blocks labelled by elements of $A$ and vertices labelled by $(T_i, \vfi: T_i \to T_j)$, where $T_j$ is the label of the block the vertex is attached to and $\vfi$ is a morphism of $G$-sets.

    Let $\cC$ be a group-theoretical subcategory of $\Rel(\cA)$. As in \Cref{ss:existence-of-abelian-envelopes}, we can again define the rank $R_f$ of an injection $f: [\mathbf{v}] \to [\mathbf{w}]$ as an object $[\mathbf{u}]$ of minimal length such that there exists a relation $[\mathbf{w}] \xhookleftarrow{f} [\mathbf{v}] \twoheadleftarrow [\mathbf{b}] \to [\mathbf{u}]$ in $\cC$. Then \Cref{a:rank-object} is satisfied if for all relations of this form the underlying partition diagram has no lower blocks. We do not know if this is true for all choices of the group-theoretical subcategory $\cC$, say, after excluding $\gtcategory{\upperdots{2}}$ as underlying group-theoretical category of partitions.
\end{remark}

\section{Groups in \texorpdfstring{$\uRep(\GL_t)$ and $\uRep(\Aff_t)$}{Rep(GLt) and Rep(Afft)}} \label{s:group-in-rep-glt}

In this section, we consider the regular category of finite dimensional vector spaces over a finite field $\FF = \FF_q$ and of finite dimensional affine spaces over $\FF$. Their skeleta, denoted by $\cA$ and $\cB$ respectively, consist of objects $\FF^n$ for $n \in \NN_0$. For each $t \in \kk$, we have a degree function $\delta_t$ sending a surjective linear or affine map $\eps: \FF^n \tto \FF^m$ to $\delta_t(\eps) = t^{n-m}$.

We first examine the relation between $\Rel(\cA)$ and $\Rel(\cB)$. We denote an affine map $\FF^n \to \FF^m$ by $f(-)+t$, where $f: \FF^n \to \FF^m$ is linear and $t \in \FF^m$ is a translation vector.

\begin{lemma} \label{l:bijection-rela-relb}
    The map
    \[
        \Phi: \Rel(\cB) \longrightarrow \Rel(\cA), \qquad \left(\FF^n \xleftarrow{f(-)+s} \FF^b \xrightarrow{g(-)+t} \FF^m\right) \longmapsto \left(\FF^n \xleftarrow{f} \FF^b \xrightarrow{g} \FF^m\right)
    \]
    induces a bijection between the sets of group-theoretical subcategories of $\Rel(\cB)$ and $\Rel(\cA)$, respectively.
\end{lemma}
\begin{proof}
    This follows from the fact that a group-theoretical subcategory $\cC \subseteq \Rel(\cB)$ must contain all relations of the form $\FF^n \xleftarrow{=} \FF^n \xrightarrow{+t} \FF^n$ for all $n \in \NN_0$ and $t \in \FF^n$. These relations are automorphisms in $\Rel(\cB)$ and we have
    \[
        \left(\FF^n \xleftarrow{f(-)+s} \FF^b \xrightarrow{g(-)+t} \FF^m\right) = 
        \left(\FF^m \xleftarrow{=} \FF^m \xrightarrow{+t} \FF^m\right) \circ
        \left(\FF^n \xleftarrow{f} \FF^b \xrightarrow{g} \FF^m\right) \circ
        \left(\FF^n \xleftarrow{+s} \FF^n \xrightarrow{=} \FF^n\right)
    \]
    for all relations $\FF^n \xleftarrow{f} \FF^b \xrightarrow{g} \FF^m$ in $\Rel(\cA)$ and all $s \in \FF^n, t \in \FF^m$.
\end{proof}

Recall from \Cref{ss:one-atom} that the relation
\[
    r \coloneqq \left(\FF^2 \xleftarrow{=} \FF^2 \to 0\right)
\]
plays a special role analogous to the partition $\upperdots{2}$ for $\uRep(S_t)$.

\begin{lemma} \label{l:single-bottom-contains-r}
    Let $\cC \subseteq \Rel(\cA)$ be a group-theoretical subcategory. If $\cC$ contains a relation of the form $\FF^n \xleftarrow{f} \FF^b \xrightarrow{g} \FF$ with $f$ not injective, then $r \in \cC$. The same holds for a group-theoretical subcategory of $\Rel(\cB)$.
\end{lemma}
\begin{proof}
    By \Cref{l:bijection-rela-relb}, it suffices to work in $\Rel(\cA)$.
    Let $0 \ne v \in \ker(f)$. Since $f, g$ are jointly injective, we have $g(v) \ne 0$. We first examine the case where $n = 2m$ is even. Denote by $c^{\o m} \in \cC(0, \FF^n)$ the $m$-fold tensor product of the cap and consider its composition with the relation $\FF^n \xleftarrow{f} \FF^b \xrightarrow{g} \FF$:
    \[\begin{tikzcd}
        && \FF^m \times_{\FF^n} \FF^b \arrow{ldd} \arrow{rdd} \arrow[two heads]{d} && \\
        && V \arrow[hook]{d} && \\
        & \FF^m \arrow{ld} \arrow[pos=.7,swap]{rd}{h} & 0 \times \FF \arrow[dotted]{lld}\arrow[dotted]{rrd} & \FF^b \arrow[pos=.7]{ld}{f} \arrow{rd}{g} & \\
        0 && \FF^n && \FF \\
    \end{tikzcd}\]
    The fiber product $\FF^m \times_{\FF^n} \FF^b$ contains the vector $(0, v)$. Its image in $V$ is the non-zero vector $(0, g(v))$. Thus the composition is the relation $0 \leftarrow \FF \xrightarrow{=} \FF$. As before, this shows that $\cC$ contains all relations.

    If $n = 2m+1$ is odd, consider the composition of $c^{\o m} \o \id_{\FF}$ with the relation $\FF^n \xleftarrow{f} \FF^b \xrightarrow{g} \FF$:
    \[\begin{tikzcd}
        && \FF^{m+1} \times_{\FF^n} \FF^b \arrow{ldd} \arrow{rdd} \arrow[two heads]{d} && \\
        && V \arrow[hook]{d} && \\
        & \FF^{m+1} \arrow{ld} \arrow[pos=.7,swap]{rd}{h} & \FF \times \FF \arrow[dotted]{lld}\arrow[dotted]{rrd} & \FF^b \arrow[pos=.7]{ld}{f} \arrow{rd}{g} & \\
        \FF && \FF^n && \FF. \\
    \end{tikzcd}\]
    Again, the fiber product $\FF^{m+1} \times_{\FF^n} \FF^b$ contains the vector $(0, v)$, so $V$ contains the non-zero vector $(0, g(v))$. Thus the composition is either the relation $\FF \leftarrow \FF^2 \rightarrow \FF$ (a rotation of $r$) or the relation $\FF \xleftarrow{0} \FF \xrightarrow{=} \FF$. In the latter case, $\cC$ must also contain $\FF^2 \xleftarrow{0} \FF^2 \xrightarrow{=} \FF^2$. Composing again with a cup, we see that $0 \leftarrow \FF^2 \xrightarrow{=} \FF^2$ is in $\cC$, as desired.
\end{proof}

\begin{lemma} \label{l:vec-fq-satisfies-assumption}
    Assume that $\FF$ contains more than two elements.
    Let $\cC \subseteq \Rel(\cA)$ or $\cC \subseteq \Rel(\cB)$ be a group-theoretical subcategory not containing $r$. Then $\cC$ satisfies \Cref{a:rank-object}.
\end{lemma}
\begin{proof}
    Again, it suffices to work in $\Rel(\cA)$: It is easy to verify that if $\cC \subseteq \Rel(\cB)$ is group-theoretical and $\Phi(\cC) \subseteq \Rel(\cA)$ satisfies \Cref{a:rank-object}, then so does $\cC$.

    We therefore fix a group-theoretical subcategory $\cC \subseteq \Rel(\cA)$ and an injective linear map $f: \FF^n \hookrightarrow \FF^m$ and define $R_f \coloneqq \FF^{r_f}$, where
    \[
        r_f \coloneqq \min\{i \in \NN_0 \mid \exists\text{relation in $\cC$ of the form } \FF^m \xhookleftarrow{f} \FF^n \twoheadleftarrow \FF^b \to \FF^i\}.
    \]
    This minimum is well-defined: The set is non-empty because the relation $\beta_f$ given by $\FF^m \xhookleftarrow{f} \FF^n \xhookrightarrow{f} \FF^m$ is a coarsening of $\id_{\FF^m}$ and therefore contained in $\cC$.

    We need to show that $R_f$ satisfies the properties (A1) and (A2) in \Cref{a:rank-object}. Consider a relation in $\cC$ of the form $\FF^m \xhookleftarrow{f} \FF^n \stackrel{\eps}{\twoheadleftarrow} \FF^b \xrightarrow{g} \FF^i$. We assume that $\eps$ is not a bijection, and deduce $i > r_f$. By \Cref{l:single-bottom-contains-r}, we have $i \geq 2$. Note that by definition, $\cC$ contains the relation $\FF \xleftarrow{\cdot a} \FF \xrightarrow{=} \FF$ for all $a\in \FF \setminus \{0\}$. It also contains the rotated version
    \[
        c_a \coloneqq \left(0 \leftarrow \FF \xrightarrow{\psmat{a\\1}} \FF^2\right).
    \]
    Choose $0 \ne u \in \ker(\eps)$. Because $f\eps$ and $g$ are jointly injective, $g(u) \ne 0$. The field $\FF$ contains more than two elements, so we can choose $j \in \{1, \dots, i-1\}$ and $a \in \FF \setminus \{0\}$ such that $g(u)_j \ne a g(u)_{j+1}$ and consider the composition of $\id_{\FF^{j-1}} \o c_a^T \o \id_{\FF^{i-j-1}}$ with $\FF^m \xhookleftarrow{f} \FF^n \stackrel{\eps}{\twoheadleftarrow} \FF^b \xrightarrow{g} \FF^i$.
    \[\begin{tikzcd}[ampersand replacement=\&,column sep=2cm]
        \&\& \FF^b \times_{\FF^i} \FF^{i-1} \arrow{ldd} \arrow{rdd} \arrow[two heads]{d} \&\& \\
        \&\& V \arrow[hook]{d} \&\& \\
        \&
            \FF^b \arrow[two heads]{ld}{\eps} \arrow[pos=.7,swap]{rdd}{g} \&
            \FF^m \times \FF^{i-2} \arrow[dotted]{lldd}\arrow[dotted]{rrdd} \&
            \FF^{i-1} \arrow[pos=.7]{ldd}{\psmat{I_{j-1}&0&0\\0&a&0\\0&1&0\\0&0&I_{i-j-1}}} \arrow{rdd}{\psmat{I_{j-1}&0&0 \\ 0&0&I_{i-j-1}}} \& \\
        \FF^n \arrow[hook,swap]{d}{f} \&\&\&\& \\
        \FF^m \&\& \FF^i \&\& \FF^{i-2}. \\
    \end{tikzcd}\]
    We show that the composed relation also has $f$ as its top, i.e.~that the image of $V \hookrightarrow \FF^m \times \FF^{i-2} \to \FF^m$ is precisely the image of $f$. For a vector $v \in \FF^b$, we need to find an element of $V$ whose image in $\FF^m$ is $f\eps(v)$. Consider the vector $v' \coloneqq v + \frac{g(v)_j-ag(v)_{j+1}}{ag(u)_{j+1}-g(u)_j} u$. We have $g(v')_j = ag(v')_{j+1}$, so
    \[
        w \coloneqq (v', g(v')_1, \dots, \widehat{g(v')_j}, \dots, g(v')_i) \in \FF^b \times \FF^{i-1}
    \]
    is contained in the fiber product $\FF^b \times_{\FF^i} \FF^{i-1}$. By construction, its image in $\FF^m$ is $f\eps(v') = f\eps(v)$. This shows that $R_f$ satisfies (A2). Together with the definition of $R_f$, this implies (A1).
\end{proof}

\begin{lemma} \label{l:quasi-regular}
    The measure on the oligomorphic group corresponding to $(\cA, \delta_t)$ or $(\cB, \delta_t)$ under \Cref{c:g,c:mu} is quasi-regular.
\end{lemma}
\begin{proof}
    For the regular category $\cA$, the group is $G \coloneqq \Aut_\FF(\Omega)$ acting oligomorphically on the vector space $\Omega \coloneqq \bigcup_{i \geq 0} \FF^i$. Its dense subgroup $\bigcup_{i \geq 0} \GL_i(\FF)$ has the same smooth finitary sets and measures, so its measures are all quasi-regular by \cite{oligotensor}*{Proposition~15.14}.

    For the regular category $\cB$, use the equivalence $\cB\simeq\cA//\FF$ sending a surjection $\lambda:E\tto\FF$ to the affine space $\lambda^{-1}(1)$ (see \cite{knop}*{\S~2, Example~4}). Since $\dim\lambda^{-1}(1)=\dim E-1$, this equivalence identifies $\delta_t$ with the degree function inherited from $\cA$. Let $\mu$ be the measure on $G$ corresponding to $\delta_t$. By the proof of \Cref{p:subgroup-of-G}(b), the group and measure associated to $\cA//\FF$ are an open subgroup $U\subseteq G$ and $\mu|_U$. Choose an open subgroup $V\subseteq G$ such that $\mu|_V$ is regular. Then $\mu|_{U\cap V}$ is regular, so $\mu|_U$ is quasi-regular.
\end{proof}

\begin{proposition} \label{prop:splitting-A-B} For each of $\cT(\cA,\delta_t)$ and $\cT(\cB,\delta_t)$, there is a global splitting object which lies in the image of the embedding of any subcategory corresponding to a group-theoretical subcategory of $\Rel(\cA)$ and $\Rel(\cB)$, respectively.
\end{proposition}

\begin{proof}
    These are special cases of \Cref{p:quasi-regular-splitting} using \Cref{l:quasi-regular}.
\end{proof}

\begin{proposition} \label{p:linear-affine-ab-env}
    Assume that $\FF$ contains more than two elements and $t \in \kk$. For a group-theoretical subcategory $\cC \subseteq \Rel(\cA)$ or $\cC \subseteq \Rel(\cB)$, the category $\cC_{\delta_t}$ admits an abelian envelope with enough projectives.
\end{proposition}
\begin{proof} This is \Cref{thm:one-atom} using \Cref{prop:splitting-A-B} and \Cref{l:vec-fq-satisfies-assumption}.
\end{proof}

\bibliographystyle{alphaurl}
\bibliography{refs}

\end{document}